# Quandle homotopy invariants of knotted surfaces

Takefumi Nosaka

**Abstract**

Given a finite quandle, we introduce a quandle homotopy invariant of knotted surfaces in the 4-sphere, modifying that of classical links. This invariant is valued in the third homotopy group of the quandle space, and is universal among the (generalized) quandle cocycle invariants. We compute the second and third homotopy groups, with respect to "regular Alexander quandles". As a corollary, any quandle cocycle invariant using the dihedral quandle of prime order is a scalar multiple of the Mochizuki 3-cocycle invariant. As another result, we determine the third quandle homology group of the dihedral quandle of odd order.

## 1   Introduction

A quandle $X$ is an algebraic system satisfying axioms that correspond to the Reidemeister moves. Given a quandle $X$, Fenn, Rourke and Sanderson [FRS] defined the rack space $BX$, in analogy to the classifying spaces of groups; Further, an invariant of framed links in $S^3$ was proposed in [FRS2, §4], which is called a quandle homotopy invariant, and is valued in the second homotopy group $\pi_2(BX)$ [see [N1] for some computations of $\pi_2(BX)$]. In addition, as a modification of the homology $H_*(BX;\mathbb{Z})$, Carter et. al [CJKLS] introduced its quandle homology denoted by $H_n^Q(X;A)$, and further quandle cocycle invariants of classical links (resp. linked surfaces) using cocycles of the cohomology $H_Q^2(X;A)$ (resp. $H_Q^3(X;A)$). For its application, the homology groups $H_*(BX;A)$ and $H_*^Q(X;A)$ of some quandles $X$ have been computed [Cla, CJKS, M1, M2, N2, NP2]. Furthermore the quandle cocycle invariants were generalized to allow the cohomology $H_Q^*(X;A)$ with local coefficients [CEGS].

As for classical links, the quandle cocycle invariants were much studied (see, e.g., [NP1, RS, N1, HN]), and are known to be derived from the quandle homotopy invariant above (see [CKS, RS]). A topological meaning of some quandle cocycle invariants was understood from the study of the homotopy group $\pi_2(BX)$ [HN].

In this paper, we introduce and study a quandle homotopy invariant of oriented linked surfaces valued in a group ring $\mathbb{Z}[\pi_3^Q(BX)]$ (Definition 2.3), modifying the above quandle homotopy invariant of links in $S^3$. Here the group $\pi_3^Q(BX)$ is defined by a quotient group of the third homotopy group $\pi_3(BX)$ (see Remark 2.2 for details). Similar to the quandle homotopy invariant of classical links, that of linked surfaces is shown to be universal among the quandle cocycle invariants with local coefficients of linked surfaces (see §2.3). It is therefore significant to estimate and determine $\pi_3^Q(BX)$; Moreover, from the study of $\pi_3^Q(BX)$, we address a problem of determining for what kinds of local coefficients the associated quandle cocycle invariants pick out completely the quandle homotopy invariant.

First, we determine the free subgroup of $\pi_3^Q(BX)$ of a finite quandle $X$, using rational homotopy theory (Theorem 3.1). We show that $\pi_3^Q(BX)$ is finitely generated, and that $\pi_3^Q(BX) \otimes \mathbb{Q}$ depends only on the number $\ell$ of "the connected components" of $X$: To be precise, $\pi_3^Q(BX) \otimes \mathbb{Q} \cong \mathbb{Q}^{\ell(\ell-1)(\ell-2)/3}$. Further, we give a topological interpretation of $\pi_3^Q(BX) \otimes \mathbb{Q}$ from a view of "link bordism groups" (Remark 4.2).

Next we deal with the torsion subgroups of $\pi_3^Q(BX)$. Since it is difficult to compute $\pi_3^Q(BX)$ in general, in this paper, we confine ourselves to regular Alexander quandles $X$. Here the regular Alexander quandle $X$ is defined to be a $\mathbb{Z}[T^{\pm}]$-module with a binary operation $(x, y) \mapsto Tx + (1-T)y$ which satisfies $(1-T)X = X$ and the minimal $e$ satisfying $(T^e - 1)X = 0$ is prime to $|X|$. We first estimate roughly $\pi_3^Q(BX)$: When $X$ is a regular Alexander quandle of odd order, we obtain an exact sequence which estimates $\pi_3^Q(BX)$ from upper bounds by the homologies $H_j(BX; \mathbb{Z})$ with $j \leq 4$ (Proposition 3.4); However, under an assumption, we succeed in determining $\pi_3^Q(BX)$ in a homological context: specifically, if $H_2^Q(X; \mathbb{Z})$ vanishes and $|X|$ is odd, then $\pi_3^Q(BX)$ is isomorphic to $H_4^Q(X; \mathbb{Z})$ (Theorem 3.5). We here remark that there are many such quandles, e.g., the Alexander quandle of the form $X = \mathbb{F}_q[T]/(T - \omega)$ where $\mathbb{F}_q$ is the finite field of order $p^{2h-1}$ and $\omega \in \mathbb{F}_q \setminus \{1, 0\}$ (Remark 3.6). In particular, we determine $\pi_3^Q(BX)$ of all Alexander quandles $X$ of prime order (see Corollary 3.7), since the author determined the homology $H_4^Q(X; \mathbb{Z})$ in [N1]. As a special case, if $X = \mathbb{Z}[T]/(p, T + 1)$ which is called the dihedral quandle, then $\pi_3^Q(BX) \cong \mathbb{Z}/p\mathbb{Z}$ and any quandle cocycle invariant of $X$ turns out to be a scalar multiple of the quandle cocycle invariant $\Psi_{\theta_p}(L) \in \mathbb{Z}[\mathbb{F}_p]$ of a 3-cocycle $\theta_p \in H_Q^3(X; \mathbb{F}_p)$ called Mochizuki 3-cocycle (Corollary 3.8). In summary, the study to compute $\pi_3^Q(BX)$ picks out all non-trivial quandle cocycle invariants of $X$.

In another direction, by a similar discussion, we also compute the second homotopy group $\pi_2(BX)^1$ with respect to regular Alexander quandles $X$. Recall that $\pi_2(BX)$ is the container of the quandle homotopy invariant of framed links. Our result on $\pi_2(BX)$ is that, if $X$ is of odd order, the following exact sequence splits (Theorem 3.9):

$$0 \longrightarrow \pi_2(BX) \longrightarrow \mathbb{Z} \oplus H_3^Q(X; \mathbb{Z}) \longrightarrow H_2^Q(X; \mathbb{Z}) \wedge_{\mathbb{Z}} H_2^Q(X; \mathbb{Z}) \longrightarrow 0.$$

In particular, $\pi_2(BX)$ is a direct summand of $\mathbb{Z} \oplus H_3^Q(X; \mathbb{Z})$. In conclusion, regarding $H_3^Q(X; \mathbb{Z})$ as a second homology with local coefficients (see Remark 4.1), Theorem 3.9 implies that the quandle homotopy invariant is completely determined as a linear sum of quandle cocycle invariants through $H_Q^2(X; A) \oplus H_Q^3(X; A)$ of the local system (Remark 5.5).

Furthermore, the above sequence enables us to compute $\pi_2(BX)$. For several quandles $X$ of order $\leq 9$, we determine exactly $\pi_2(BX)$ (see Table 1), following the values of $H_2^Q(X; \mathbb{Z})$ and $H_3^Q(X; \mathbb{Z})$ presented in [LN]. Furthermore, since Mochizuki [M1, M2] has computed $H_Q^2(X; \mathbb{F}_q) \oplus H_Q^3(X; \mathbb{F}_q)$ for $X = \mathbb{F}_q[T]/(T - \omega)$, we determine $\pi_2(BX) \otimes_{\mathbb{Z}} \mathbb{F}_p$ and have made the quandle homotopy invariants computable (see §A). One of the noteworthy results is that if $X$ is the product $h$-copies of the dihedral quandle, i.e., $X = \left(\mathbb{Z}[T]/(p, T + 1)\right)^h$, then the dimension of $\pi_2(BX) \otimes_{\mathbb{Z}} \mathbb{F}_p$ is a quadruple function of $h$: to be precise, $\dim_{\mathbb{F}_p}(\pi_2(BX) \otimes_{\mathbb{Z}} \mathbb{F}_p) = 1 + (h^2(h^2+11)/12)$ (Corollary A.4).

---
[1] $\pi_2(BX) \cong \mathbb{Z} \oplus \pi_2(B^Q X)$ is known [N1], where $B^Q X$ is the quandle space. In §3 and 4, we mainly deal with $\pi_2(B^Q X)$.

As a corollary of our work, for an odd $m$, we determine $H_3^Q(X;\mathbb{Z})$ of the dihedral quandle $X$ of order $m$, i.e., $X = \mathbb{Z}[T]/(m, T+1)$. Since Hatakenaka and the author [HN] computed $\pi_2(BX) \cong \mathbb{Z} \oplus \mathbb{Z}/m\mathbb{Z}$, it follows from Theorem 3.9 above that $H_3^Q(X;\mathbb{Z}) \cong \mathbb{Z}/m\mathbb{Z}$ (Corollary 3.11). When $m$ is an odd prime, $H_3^Q(X;\mathbb{Z}) \cong \mathbb{Z}/m\mathbb{Z}$ was conjectured by Fenn, Rourke and Sanderson (see [Oht, Conjecture 5.12]), and solved by Niebrzydowski and Przytycki using homological algebra [NP1]. On the other hand, our proof is another approach from $\pi_2(BX)$, and further gives a generalization of the conjecture (see also Remark 3.12).

This paper is organized as follows. In §2, we review quandles and the rack spaces $BX$, and define the quandle homotopy invariant. In §3, we state our theorems and corollaries. In §4, we review the quandle homology groups, and prove Theorem 3.1 on the rational homotopy group $\pi_3^Q(BX) \otimes \mathbb{Q}$. In §5, we show Theorems 3.5 and 3.9. In §6, we compute $\pi_2(BX) \otimes \mathbb{Z}/p\mathbb{Z}$ of some Alexander quandles $X$ on finite fields.

**Notational convention** Throughout this paper, $p$ is a prime and $\mathbb{F}_q$ is a finite field such that $|\mathbb{F}_q| = q = p^h$. We denote by $\mathbb{Z}_m$ the cyclic group of order $m \in \mathbb{Z}$. A symbol "$\otimes$" is always the tensor product over $\mathbb{Z}$. For a group $G$, we denote by $\mathbb{Z}[G]$ the group ring of $G$ over $\mathbb{Z}$. For a $\mathbb{Z}$-module $N$, $N_{(p)}$ means the localization at a prime $p$. Further, $\Lambda^*(N)$ denotes the $*$-part of the exterior algebra over $N$. A symbol "pt." stands for a single point.

All embeddings of surfaces into $S^4$ are assumed to be $C^\infty$-class. Furthermore, we assume that a surface is connected, oriented and closed.

## 2 Quandle homotopy invariants of linked surfaces

We review quandles and the rack spaces in §2.1, and introduce quandle homotopy invariants of linked surfaces in §2.2.

### 2.1 Review of quandles and of rack spaces

A *quandle* is a set $X$ with a binary operation $(x,y) \to x * y$ such that, for any $x, y, z \in X$, $x * x = x$, $(x * y) * z = (x * z) * (y * z)$ and there exists uniquely $w \in X$ satisfying $w * y = x$. For example, a $\mathbb{Z}[T^{\pm}]$-module $M$ has a quandle structure given by $x * y := Tx + (1-T)y$, called an *Alexander quandle*. For a quandle $X$, *the associated group* of $X$ is defined by the group presentation $\mathrm{As}(X) := \langle x \in X \mid x \cdot y = y \cdot (x * y) \rangle$. Remark that, for any quandle $X$, $\mathrm{As}(X)$ is of infinite order, since $\mathrm{As}(X)$ has an epimorphism onto $\mathbb{Z}$ obtained by the length of words. An *$X$-set* is defined to be a set $Y$ with an action of $\mathrm{As}(X)$. For example, when $Y = X$, the quandle $X$ is an $X$-set itself by the quandle operation.

We next review $X$-colorings [CJKLS, Definition 5.1]. Let $D$ be a broken diagram of a linked surface $L$. In this paper, a broken diagram is in $S^3$. An *$X$-coloring* is a map $C$ from the set of sheets of $D$ to $X$ such that, for every double-point curve as shown in the left of Figure 1, the three sheets $\alpha, \beta, \gamma$ satisfy $C(\gamma) = C(\alpha) * C(\beta)$. This definition is compatible with triple-points (see the right of Figure 1). We denote the set consisting of such $X$-colorings of $D$ by $\mathrm{Col}_X(D)$. It is widely known (see, e.g., [CJKLS, Theorem 5.6]) that if two diagrams $D$ and $D'$ are related by the Roseman moves, then there naturally exists a 1-1 correspondence between $\mathrm{Col}_X(D)$ and $\mathrm{Col}_X(D')$. Furthermore, when $X$ is an Alexander quandle of the form $X = \mathbb{Z}_p[T^{\pm}]/(h(T))$ and $L$ is a knotted surface, the set $\mathrm{Col}_X(D)$ is obtained from the Alexander polynomials of $L$ (see [Ino] for details).

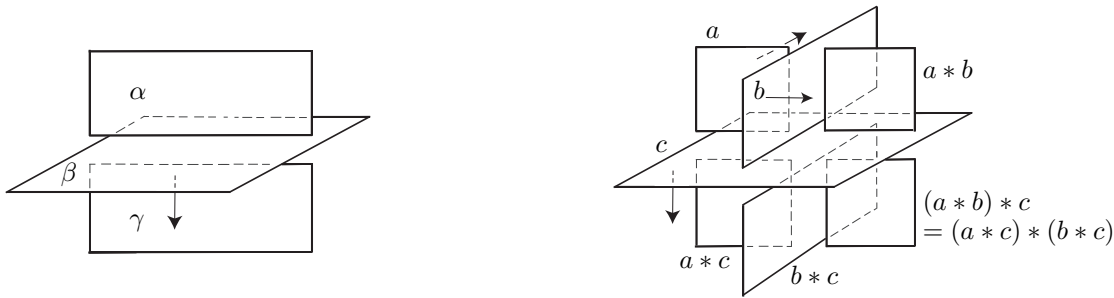

Figure 1: Coloring conditions around double-point curves and triple-points. Here $a, b, c \in X$

Finally, we review the rack space introduced by Fenn-Rourke-Sanderson [FRS]. Fix an $X$-set $Y$. Equipping $X$ and $Y$ with their discrete topology, we start with $\bigcup_{n \geq 0}(Y \times ([0,1] \times X)^n)$, and consider the equivalence relations given by

$$(y; t_1, x_1, \ldots, x_{j-1}, 1, x_j, t_{j+1}, \ldots, t_n, x_n) \sim (y \cdot x_j; t_1, x_1 * x_j, \ldots, t_{j-1}, x_{j-1} * x_j, t_{j+1}, x_{j+1}, \ldots, t_n, x_n),$$

$$(y; t_1, x_1, \ldots, x_{j-1}, 0, x_j, t_{j+1}, \ldots, t_n, x_n) \sim (y; t_1, x_1, \ldots, t_{j-1}, x_{j-1}, t_{j-1}, x_{j+1}, \ldots, t_n, x_n).$$

Here our description divides the element of $Y$ from the elements of $([0,1] \times X)^n$ by a semicolon ';'. Then the *rack space* $B_Y X$ is defined by the quotient space (we later explain the 3-skeleton for more details in §2.2). By construction, we have a cell decomposition of $B_Y X$, by regarding the projection $\bigcup_{n \geq 0}(Y \times ([0,1] \times X)^n) \to B_Y X$ as the characteristic maps. We easily see that if the action on $Y$ of $\mathrm{As}(X)$ is transitive, the space $B_Y X$ is path-connected (see [Cla, Proposition 7] for details). When $Y$ is a single point, we denote the path-connected space $B_Y X$ by $BX$ for short. We remark that $\pi_1(BX) \cong \mathrm{As}(X)$ by the definition of the 2-skeleton of $BX$.

## 2.2 Definition of a quandle homotopy invariant of linked surfaces

Our purpose is to define a quandle homotopy invariant of linked surfaces. This is an analogue of the quandle homotopy invariant of (framed) links in $S^3$ [FRS2, §3] (see also [N1, §2]). For the purpose, given an $X$-set $Y$, we first define a quandle space as follows. Define a subspace by

$$\bigcup_{n \geq 2} \{(y; t_1, x_1, \ldots, t_n, x_n) \in Y \times ([0,1] \times X)^n \mid x_i = x_{i+1} \text{ for some } i \in \mathbb{Z}, \ i \leq n-1\}. \quad (1)$$

We denote this by $X_Y^D$, and consider a composite $\iota_Y^D : X_Y^D \hookrightarrow \bigcup_{n \geq 0}(Y \times ([0,1] \times X)^n) \to B_Y X$. Here the second map is the projection mentioned above.

**Definition 2.1.** Let $X$ be a quandle. Let $Y$ be an $X$-set. We define a space $B_Y^Q X$ by the cone of the composite map $\iota_Y^D : X_Y^D \to B_Y X$. We refer to the space as *quandle space*.

When $Y$ is a single point, we denote $B_Y^Q X$ by $B^Q X$ for short. The homology $H_*(B^Q X; \mathbb{Z})$ coincides with the homology of the quandle complex introduced in [CJKLS] (see also §4.1). Roughly speaking, the space $B^Q X$ is a geometric realization of the quandle complex.

We now explain the 3-skeleton of the rack space $BX$ and the 4-skeleton of the quandle space $B^Q X$ in more details. The 1-skeleton of $BX$ is a bouquet of $|X|$-circles labeled by elements of $X$. The 2-skeleton of $BX$ is obtained from the 1-skeleton by attaching 2-cells of squares labeled by $(a,b)$ for any $a,b \in X$, where the 4 edges with $X$-labels as shown in Figure 2 are attached to the corresponding 1-cells. In addition, we attach $|X^3|$-cubes labeled by $(a,b,c) \in X^3$, whose six faces are labeled as shown in Figure 2, to the corresponding 2-cells. The resulting space is the 3-skeleton of $BX$.

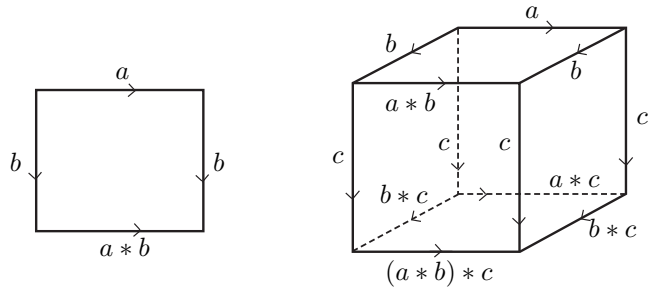

Figure 2: The 2-cell labeled by $(a,b)$ and the 3-cell labeled by $(a,b,c)$

We next explain the 4-skeleton of $B^Q X$. To begin, the 3-skeleton of $B^Q X$ is obtained from the 3-skeleton of $BX$ by attaching 3-cells which bound 2-cells labeled by $(a,a)$ for $a \in X$. The 4-skeleton of $B^Q X$ is obtained from the 3-skeleton of $B^Q X$ by attaching 4-cells as follows: for $a, b \in X$, a 3-cell labeled by either $(a,a,b)$ or $(a,b,b)$ is bounded by a 4-cell, and the 4-dimensional cube labeled by $(a,b,c,d) \in X^4$ bounds the eight corresponding 3-cells as shown in Figure 3. Remark that the precceding quandle space considered in [N1, §2] [2] is the 3-skeleton of the space $B^Q X$, and the higher skeleton was not defined in [N1].

In order to construct an invariant of linked surfaces, we prepare a cell decomposition of $S^3$. Let $L \subset S^4$ be a linked surface. We fix a broken surface diagram $D \subset S^3$ of $L$. Regarding

---
[2] In [N1], the rack (resp. quandle) space is denoted by $\widehat{BX}$ (resp. $BX$).

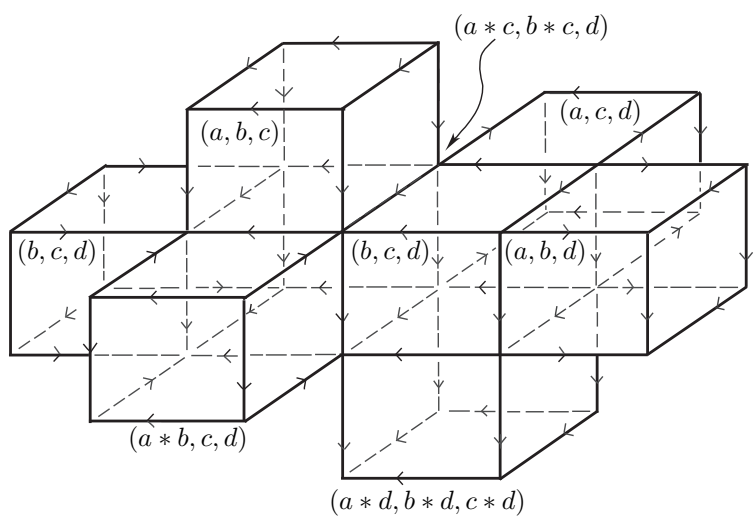

Figure 3: The eight cubes with labels as the boundary of the 4-cell labeled by $(a, b, c, d)$

$D$ as a decomposition of $S^3$ by an immersed surface, we consider the dual decomposition. Recall that the broken diagram is locally composed of some double-point curves, branch-points and triple-points (see, e.g., [CS] [CJKLS, §5] for the definitions). We then modify the dual decomposition around each branch-points $P_B$ of $D$ as follows. To begin, we add an interval which connects between $y$ and $z$ shown in the middle of Figure 4. Put a square $\mathcal{D}_{xyz}$ which bounds the three intervals, and further attach a 3-cell which bounds the square $\mathcal{D}_{xyz}$ and contains the branch-point $P_B$ (see the right of Figure 4); we here notice that the 3-cell forms a cone on a square. In addition, we attach a cube to the square $\mathcal{D}_{xyz}$ from the opposite direction. For example, if the Whitney umbrella crosses another sheet such as the left of Figure 5, we here attach 2-cells and 3-cells to the dual decomposition shown as the right of Figure 5. We orient the modified cell decomposition by using the orientations of $D$.

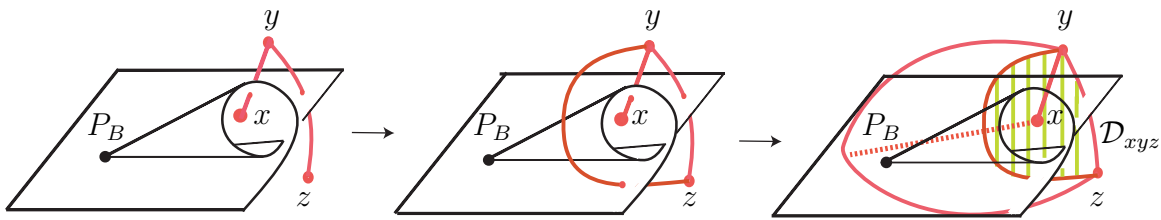

Figure 4: The modification of the dual decomposition I

Next, given an $X$-coloring $C$ of $D$, we will construct a map $\xi_{D,C} : S^3 \to B^Q X$, as follows. We take the 0-cells of the modified decomposition to the single 0-cell of $BX$. We send the 1-cells at sheets colored by $a$ to the 1-cell labeled by $a$. We take the 2-cells around the double-point curves colored by $(a,b)$ to the 2-cell labeled by $(a,b)$. We map the 3-cells around the triple-points colored by $(a,b,c)$ to the 3-cell labeled by $(a,b,c)$ as shown in Figure 6. Further, we take the modified 3-cells around the branch-points colored by $a \in X$ to the above cone on the 2-cell labeled by $(a,a)$ (see the right of Figure 4). By collecting them, we obtain a cellular map $\xi_{D,C} : S^3 \to B^Q X$. We denote the homotopy class of $\xi_{D,C}$ by $\Xi_X(D;C) \in \pi_3(B^Q X)$. By the construction of 4-cells of $B^Q X$, we can verify that the homotopy class $\Xi_X(D;C)$ is invariant under the Roseman moves (cf. [CJKLS, Theorem 5.6]). For example, the invariance of the move in the bottom right of [CJKLS, Figure 7] corresponds to the 4-cell in Figure 3.

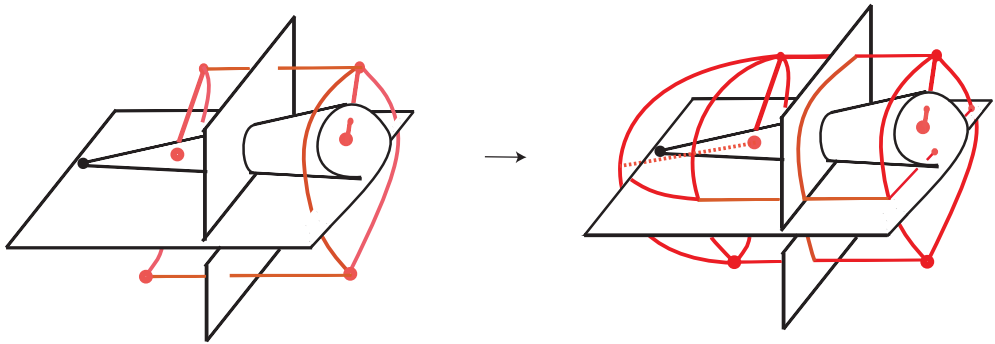

Figure 5: The modification of the dual decomposition II

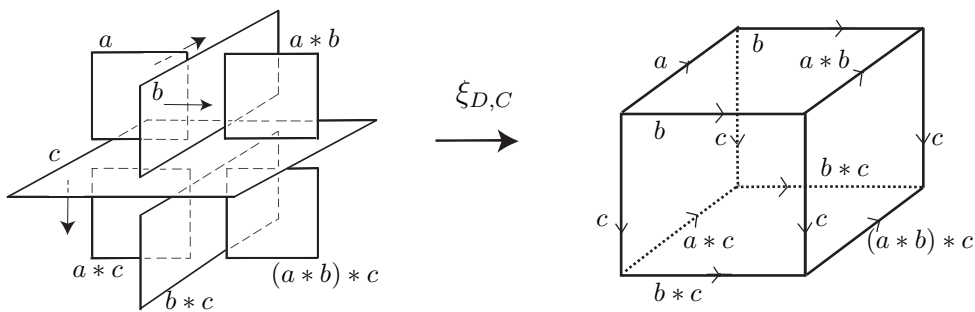

Figure 6: From a triple-point colored by $(a,b,c)$ to the 3-cell labeled by $(a,b,c)$

**Remark 2.2.** It is known that any oriented surface-link is represented by a broken diagram $D$ without branch-point (see [CS]). Recall that we defined the group $\pi_3^Q(BX)$ to be the quotient of $\pi_3(BX)$ by a normal subgroup that takes account of branch-points of knotted surfaces. However, for such a diagram $D$, we claim that the map $\xi_{D,C} : S^3 \to B^Q X$ misses the cell of $X_Y^p$. Actually, the map $\xi_{D,C}$ can, by definition, be constructed by only the usual dual composition of $S^3$ by $D$. Hence, the map $\xi_{D,C} : S^3 \to B^Q X$ factors through the rack space $BX$, which implies that the homotopy class $\Xi_X(D;C)$ is derived from $\pi_3(BX)$. Therefore, for the study of the invariant, it suffices to analyze the image of $(i_X)_* : \pi_3(BX) \to \pi_3(B^Q X)$, where $i_X$ is the inclusion $BX \to B^Q X$. Let $\pi_3^Q(BX)$ denote the image $\mathrm{Im}((i_X)_*)$.

We then define

**Definition 2.3.** Let $X$ be a finite quandle. We define the group $\pi_3^Q(BX)$ by the image of the induced map $\pi_3(BX) \to \pi_3(B^Q X)$ by the inclusion $BX \to B^Q X$. Further we define a *quandle homotopy invariant* of a linked surface $L$ by the following expression:
$$\Xi_X(L) \ = \ \sum_{C \in \mathrm{Col}_X(D)} \Xi_X(D;C) \ \in \mathbb{Z}[\pi_3^Q(BX)].$$
Here $D$ is a broken diagram of $L$, and $\mathrm{Col}_X(D)$ denotes the set of $X$-colorings of $D$.

In §3.1, we study the container $\pi_3^Q(BX)$.

## 2.3 Reconstucion of generalized quandle cocycle invariants

We will reformulate the (generalized) quandle cocycle invariants in [CJKLS, CKS, CEGS] from our quandle homotopy invariant. We first remark that, if $B_Y X$ is path-connected, the inclusion $B_Y X \hookrightarrow B_Y^Q X$ induces $\pi_1(B_Y X) \cong \pi_1(B_Y^Q X)$, since $X_Y^D$ given in (1) has no 1-cell by definition. For an action of $\pi_1(B^Q X) = \mathrm{As}(X)$ on an abelian group $A$, we put a 3-cocycle $\psi \in H^3(B^Q X; A)$ with local coefficients. The homology $H_3(B^Q X; A)$ coincides with the quandle homology $H_3^Q(X; A)$ introduced in [CEGS, §2 and §7] (see also §4.1). Let $\mathfrak{H} : \pi_3(B^Q X) \to H_3(B^Q X; A)$ be the Hurewicz homomorphism with local coefficients. Then the *quandle cocycle invariant* of $L$ is defined by

$$\Phi_\psi(L) = \sum_{C \in \mathrm{Col}_X(D)} \langle \psi, \mathfrak{H}(\Xi_X(D;C)) \rangle \in \mathbb{Z}[A]. \tag{2}$$

This equality implies that the quandle homotopy invariant is universal among the quandle cocycle invariants.

We give some remarks on the formula (2). It is not difficult to see that the formula (2) coincides with the combinatorial definition of the (generalized) quandle cocycle invariants in [CEGS, §7], by the definition of the Hurewicz homomorphism and the construction of $\xi_{D,C}$. Hence, the equality (2) means that, if knowing a concrete presentation of $\psi$, we can calculate some parts of the quandle homotopy invariant. However there are many choices of actions $\pi_1(BX) \curvearrowright A$; then we will approach a problem of determining for what kinds of local systems the associated quandle cocycle invariants pick out the quandle homotopy invariant. §3.1 gives an answer for some Alexander quandles.

Finally, we discuss a local system. Regarding $X$ itself as an $X$-set, the free $A$-module $A\langle X \rangle$ generated by $x \in X$ is a $\pi_1(BX)$-module. We then consider a 4-cocycle $\psi \in H^4(B^Q X; A)$ to be a cocycle of $H^3(B^Q X; A\langle X \rangle)$ (see Remark 4.1 for detail); The quandle cocycle invariant of $\psi$ coincides with the *shadow cocycle invariant* described in [CKS, §5].

## 2.4 Properties of quandle homotopy invariants

We state some properties of quandle homotopy invariants of linked surfaces, similar to those of classical links in [N1, §5]. For this, we review connected components of a quandle $X$. By definition of $X$, we note that $(\bullet * y) : X \to X$ is bijective for any $y \in X$. Let $\mathrm{Inn}(X)$ denote the subgroup of $\mathfrak{S}_{|X|}$ generated by the all right actions $(\bullet * y)$, and is called the *inner automorphism group* of $X$. The *connected components* of a quandle $X$ are defined by the orbits of the action of $\mathrm{Inn}(X)$ on $X$. For example, we let $T_\ell$ be composed $\ell$-points with a binary operation defined by $x * y = x$ for $x, y \in T_\ell$; then the number of the connected components of $T_\ell$ is $\ell$. The quandle $T_\ell$ is called a *trivial quandle*. Furthermore, a quandle $X$ is said to be *connected*, if the right action of $\mathrm{Inn}(X)$ on $X$ is transitive. For example, it is known [LN, Proposition 1] that a finite Alexander quandle is connected if and only if $(1-T)$ is invertible.

We give a slight reduction of quandle homotopy invariants:

**Proposition 2.4.** *(cf. [N1, Lemma 5.5]) Let $X$ be a finite connected quandle, $x$ an element of $X$ and $D$ a link diagram of a linked surface $L$. We fix a sheet of $D$, and denote $\mathrm{Col}_X^x(D)$ by the subset of $\mathrm{Col}_X(D)$ satisfying that the sheet is colored by $x \in X$. Then,*

$$\Xi_X(L) = |X| \sum_{C \in \mathrm{Col}_X^x(D)} \Xi_X(D; C) \in \mathbb{Z}[\pi_3^Q(BX)]. \qquad (3)$$

*Proof.* For any $y \in X$, We see that the bijection $(\bullet * y) : X \to X$ induces a bijection $T_y : \mathrm{Col}_X^x(D) \cong \mathrm{Col}_X^{x*y}(D)$ satisfying that $\Xi_X(D; C) = \Xi_X(D; T_y(C)) \in \pi_3^Q(BX)$ (see the proof of [FRS2, Proposition 5.2]). Therefore, since $X$ is connected, for any $z \in X$, we have

$$\sum_{C \in \mathrm{Col}_X^x(D)} \Xi_X(D; C) = \sum_{C \in \mathrm{Col}_X^z(D)} \Xi_X(D; C) \in \mathbb{Z}[\pi_3^Q(BX)],$$

which immediately results the required equality (3). □

Furthermore, we can show the formulas for the connected sum and the mirror image of the quandle homotopy invariant. We give the formulas without proofs, since the proofs are similar to the discussions in [N1, §5] on the quandle homotopy invariant of classical links.

**Proposition 2.5.** *(cf. [N1, Proposition 5.1]) Let $X$ be a connected quandle of finite order. Then, for knotted surfaces $K_1$ and $K_2$,*

$$\Xi_X(K_1 \# K_2) = \frac{1}{|X|} \Xi_X(K_1) \cdot \Xi_X(K_2) \ \in \mathbb{Z}[\pi_3^Q(BX)],$$

*where $K_1 \# K_2$ denotes the connected sum of $K_1$ and $K_2$.*

**Proposition 2.6.** *(cf. [N1, Proposition 5.7]) Let $L$ be a linked surface. Let $-L^*$ denote the mirror image of $L$ with opposite orientation. For any finite quandle $X$, $\Xi_X(-L^*) = \iota\big(\Xi_X(L)\big)$, where $\iota$ is the map of $\mathbb{Z}[\pi_3^Q(BX)]$ induced by $\iota(x) = x^{-1}$ for any $x \in \pi_3^Q(BX)$.*

The similar formulas to Propositions 2.4, 2.5 hold for the quandle cocycle invariants via the equality (2).

# 3 Results of $\pi_3^Q(BX)$ and of $\pi_2(B^Q X)$

In §2.2, we defined an invariant valued in $\mathbb{Z}[\pi_3^Q(BX)]$. We now state our results on the free and torsion parts of the group $\pi_3^Q(BX)$ in §3.1. We also compute $\pi_2(B^Q X)$ in §3.2.

## 3.1 Results about $\pi_3^Q(BX)$

To begin with, we determine the free subgroup of $\pi_3^Q(BX)$ as follows:

**Theorem 3.1.** *Let $X$ be a finite quandle with $\ell$-connected components (see §2.4 for the definition). Then $\pi_3^Q(BX)$ is finitely generated. Further, the rational homotopy groups are*

$$\dim_{\mathbb{Q}}\bigl(\pi_2(B^Q X) \otimes \mathbb{Q}\bigr) = \frac{\ell^2 - \ell}{2}, \qquad \dim_{\mathbb{Q}}\bigl(\pi_3^Q(BX) \otimes \mathbb{Q}\bigr) = \frac{\ell(\ell-1)(\ell-2)}{3}.$$

*In particular, if $X$ is connected, $\pi_3^Q(BX)$ is a finite abelian group.*

We later discuss a topological meaning of the quandle homotopy invariant after tensoring with $\mathbb{Q}$ valued in $\pi_3^Q(BX) \otimes \mathbb{Q} \cong \mathbb{Q}^{\frac{\ell(\ell-1)(\ell-2)}{3}}$ (Remark 4.2).

Furthermore, we give a corollary:

**Corollary 3.2.** *Let $X$ be a finite connected quandle. Then, for any 3-cocycle $\phi \in H_Q^3(X;\mathbb{Z})$ with local coefficients, the quandle cocycle invariant $\Phi_\phi(L)$ is trivial. Namely, $\Phi_\phi(L) \in \mathbb{Z}[0_{\mathbb{Z}}]$.*

*Proof.* Notice $\pi_3^Q(BX) \otimes \mathbb{Q} \cong 0$ by Theorem 3.1. Since $\Phi_\phi(L)$ is derived from the ring $\mathbb{Z}[\pi_3^Q(BX) \otimes \mathbb{Q}]$ by (2), the value of $\Phi_\phi(L)$ is trivial. □

**Remark 3.3.** In [CEGS, Tables 3,4 and 5], for a connected quandle of order 3, a non-trivial quandle cocycle invariant $\Phi_\kappa(L)$ valued in $\mathbb{Z}[\mathbb{Z}^3]$ was proposed. These values are, however, incorrect, e.g., compare Propositions 2.4 and 2.6 with the tables.

Our next step is to study the torsion subgroup of $\pi_3^Q(BX)$. However, in general, it is difficult to compute explicitly homotopy groups of spaces. In this paper, we then confine ourselves to $\pi_3^Q(BX)$ of regular Alexander quandles $X$. Here an Alexander quandle $X$ is said to be *regular*, if $X$ is connected and the minimal number $e$ satisfying $(T^e - 1)X = 0$ is prime to $|X|$. For example, given $\omega \in \mathbb{F}_q$ with $\omega \neq 0, 1$, the Alexander quandle on $\mathbb{F}_q$ obtained by $x * y = \omega x + (1-\omega)y$ is regular, since $\omega^{q-1} = 1 \in \mathbb{F}_q$. Also, for an odd number $m \in \mathbb{Z}$, the Alexander quandle of the form $D_m := \mathbb{Z}_m[T]/(T+1)$ is regular ($D_m$ is called a *dihedral quandle*). Regular Alexander quandles were partly dealt with in [Cla, §4.1].

First, we give an estimate of $\pi_3(BX)$ by homologies of the rack space $BX$ in homological contexts. The homology of $BX$ is usually called *rack homology of $X$* (see §4.1 for details).

**Proposition 3.4.** *Let $X$ be a regular Alexander quandle of odd order. Let $H_j^R(X)$ denote the integral homology $H_j(BX;\mathbb{Z})$. Let $Y_2$ be a product of Eilenberg-MacLane Spaces $K(\pi_2(BX), 2) \times K(H_2^R(X), 1)$. Then there exists an exact sequence*

$$H_4(Y_2; \mathbb{Z}) \longrightarrow \pi_3(BX) \longrightarrow H_4^R(X) \longrightarrow H_3(H_2^R(X); \mathbb{Z}) \oplus (H_2^R(X) \otimes \pi_2(BX)) \longrightarrow 0.$$

*Here $H_3(H_2^R(X);\mathbb{Z})$ is the third group homology of the abelian group $H_2^R(X)$.*

Proposition 3.4 is proven by a routine calculation of Postnikov tower in §5. We see later that $\pi_2(BX)$ is determined by $H_2^R(X)$ and $H_3^R(X)$ [Theorem 3.9 and (5)]; in conclusion, $\pi_3(BX)$ and its quotient $\pi_3^Q(BX)$ are estimated from upper bounds using the rack homologies $H_j^R(X)$ with $j \leq 4$, although the exact sequence is not useful to compute $\pi_3(BX)$ and $\pi_3^Q(BX)$.

However, under an assumption, we determine $\pi_3^Q(BX)$ of regular Alexander quandles $X$:

**Theorem 3.5.** *Let $X$ be a regular Alexander quandle of odd order. Let $H_n^Q(X;\mathbb{Z})$ be the quandle homology with trivial coefficients. Assume that the second homology $H_2^Q(X;\mathbb{Z})$ vanishes. Then $\pi_3(BX) \cong H_4^Q(X;\mathbb{Z}) \oplus \mathbb{Z}/2\mathbb{Z}$ and $\pi_3^Q(BX) \cong H_4^Q(X;\mathbb{Z})$.*

**Remark 3.6.** There are many Alexander quandles satisfying the assumption. For example, we later show (Lemma A.6) that, if $X$ is of the form $X = \mathbb{F}_q[T]/(T - \omega)$ and $\mathbb{F}_q$ is an extension of odd degree, then $H_2^Q(X;\mathbb{Z}) \cong 0$.

We defer the proof until §4.2. As a corollary, we now determine $\pi_3^Q(BX)$ for all Alexander quandles $X$ of prime order:

**Corollary 3.7.** *For an odd prime $p$, we let $X$ be the Alexander quandle of the form $\mathbb{Z}[T]/(p, T - \omega)$. If $\omega \neq -1$, then $\pi_3^Q(BX) \cong 0$. If $\omega = -1$, then $\pi_3^Q(BX) \cong \mathbb{Z}_p$.*

*Proof.* Mochizuki showed $H_2^Q(X;\mathbb{Z}) \cong 0$ [M1, Corollary 2.2]. It is shown [N2, Corollary 2.3] that if $\omega \neq -1$, then $H_4^Q(X;\mathbb{Z}) \cong 0$, and that if $\omega = -1$, then $H_4^Q(X;\mathbb{Z}) \cong \mathbb{Z}_p$. □

Furthermore, we consider the dihedral quandle, i.e., the case of $\omega = -1$. The third cohomology $H_Q^3(X;\mathbb{F}_p) \cong \mathbb{F}_p$ is shown [M1, M2], and the generator $\theta_p \in H_Q^3(X;\mathbb{F}_p)$ is called *Mochizuki 3-cocycle*. Any quandle cocycle invariant using the dihedral quandle is summarized to assess the quandle cocycle invariant of the Mochizuki 3-cocycle as follows.

**Corollary 3.8.** *Let $X = \mathbb{Z}_p[T]/(T + 1)$ be the dihedral quandle of order $p$. The group $\pi_3^Q(BX) \cong \mathbb{Z}_p$ is generated by $\Xi_X(K;C)$, where $C$ is an $X$-coloring of the 2-twist spun $(2,p)$-torus knot $K$. Further, any quandle cocycle invariant is equal to a scalar multiple of the quandle cocycle invariant $\Phi_{\theta_p}(L) \in \mathbb{Z}[\mathbb{F}_p]$ using the Mochizuki 3-cocycle $\theta_p \in H_Q^3(X;\mathbb{F}_p)$.*

*Proof.* It is known [AS, Theorem 6.3] that the quandle cocycle invariant $\Phi_{\theta_p}(K)$ of $K$ is non-trivial. Hence, $\pi_3^Q(BX) \cong \mathbb{Z}_p$ is generated by $\Xi_X(K;C)$ for some $X$-coloring $C$, and the map $\langle \theta_p, \mathfrak{H}(\bullet) \rangle : \pi_3^Q(BX) \to \mathbb{F}_p$ is isomorphic. Hence, for any cocycle $\psi \in H_Q^3(X;A)$, the quandle cocycle invariant is derived from $\Phi_{\theta_p}(L)$ via the inverse map $\langle \theta_p, \mathfrak{H}(\bullet) \rangle^{-1}$ by the formula (2). □

### 3.2 Computation of the second homotopy groups of quandle spaces $\pi_2(B^Q X)$

Changing the subject, we discuss the quandle homotopy invariant of 1-dimensional links considered in [N1]. This invariant is valued in the group ring $\mathbb{Z}[\pi_2(B^Q X)]$, and is universal among the generalized cocycle invariants in [CEGS, §6], similar to the equality (2) (see [N1, §2] for details). For a finite connected quandle $X$, it is shown [N1, Theorem 3.6 and Proposition 3.12] that $\pi_2(B^Q X)$ is finite and $\pi_2(BX) = \pi_2(B^Q X) \oplus \mathbb{Z}$; The author further computed $\pi_2(B^Q X)$ for connected quandles of order $\leq 6$ [N1, §5].

In this paper, we explicitly compute $\pi_2(B^Q X)$ of regular Alexander quandles $X$ of odd order as the following theorem which we will prove in §5s:

**Theorem 3.9.** *Let $X$ be a regular Alexander quandle of finite order.*

(i) If $H_2^Q(X;\mathbb{Z})$ vanishes, then $\pi_2(B^Q X)$ is isomorphic to $H_3^Q(X;\mathbb{Z})$.

(ii) If $|X|$ is odd, the following exact sequence splits:

$$0 \longrightarrow \pi_2(B^Q X) \longrightarrow H_3^Q(X;\mathbb{Z}) \longrightarrow \Lambda^2\big(H_2^Q(X;\mathbb{Z})\big) \longrightarrow 0. \tag{4}$$

Hence we give a stronger estimate of the torsion part of $\pi_2(B^Q X)$ than [N1, Theorem 3.6].

**Corollary 3.10.** *Let $X$ be a regular Alexander quandle of odd order. Any element of $\pi_2(B^Q X)$ is annihilated by $|X|$.*

*Proof.* It is known [N2, Corollary 6.2] that $H_3^Q(X;\mathbb{Z})$ is annihilated by $|X|$. □

We next give two applications. First, we explicitly compute homotopy groups $\pi_2(B^Q X)$ for some regular Alexander quandles $X$. If we concretely know $H_2^Q(X;\mathbb{Z})$ and $H_3^Q(X;\mathbb{Z})$, then Theorem 3.9 enables us to compute $\pi_2(B^Q X)$. For example, the computations of all regular Alexander quandles $X$ with $5 \leq |X| \leq 9$ are listed in Table 1 below. In Table 1, the results of $H_2^Q(X;\mathbb{Z})$ and of $H_3^Q(X;\mathbb{Z})$ follow from [LN, Table 1]. For more example, in §A we present a computation of $\pi_2(B^Q X) \otimes \mathbb{Z}_p$ for all Alexander quandles on $\mathbb{F}_q$ with $p > 2$.

| $X$ | $H_2^Q(X;\mathbb{Z})$ | $H_3^Q(X;\mathbb{Z})$ | $\pi_2(B^Q X)$ |
|---|---|---|---|
| $\mathbb{Z}[T]/(5, T+1)$ | 0 | $\mathbb{Z}/5\mathbb{Z}$ | $\mathbb{Z}/5\mathbb{Z}$ |
| $\mathbb{Z}[T]/(5, T-\omega)$ | 0 | 0 | 0 |
| $\mathbb{Z}[T]/(7, T+1)$ | 0 | $\mathbb{Z}/7\mathbb{Z}$ | $\mathbb{Z}/7\mathbb{Z}$ |
| $\mathbb{Z}[T]/(7, T-\omega)$ | 0 | 0 | 0 |
| $\mathbb{Z}[T]/(2, T^3+T^2+1)$ | 0 | $\mathbb{Z}/2\mathbb{Z}$ | $\mathbb{Z}/2\mathbb{Z}$ |
| $\mathbb{Z}[T]/(2, T^3+T+1)$ | 0 | $\mathbb{Z}/2\mathbb{Z}$ | $\mathbb{Z}/2\mathbb{Z}$ |
| $\mathbb{Z}[T]/(9, T+1)$ | 0 | $\mathbb{Z}/9\mathbb{Z}$ | $\mathbb{Z}/9\mathbb{Z}$ |
| $\mathbb{Z}[T]/(3, T^2+1)$ | $\mathbb{Z}/3\mathbb{Z}$ | $(\mathbb{Z}/3\mathbb{Z})^3$ | $(\mathbb{Z}/3\mathbb{Z})^3$ |
| $\mathbb{Z}[T]/(3, T^2+T-1)$ | 0 | 0 | 0 |
| $\mathbb{Z}[T]/(3, T^2-T-1)$ | 0 | 0 | 0 |

Table 1: Some homotopy groups $\pi_2(B^Q X)$ obtained from $H_2^Q(X;\mathbb{Z})$ and $H_3^Q(X;\mathbb{Z})$. Here $\omega \neq \pm 1$.

In another application, the result of $\pi_2(BX)$ determines a third quandle homology:

**Corollary 3.11.** *(cf. [NP1, §3]) For an odd $m$, let $D_m = \mathbb{Z}[T]/(m, T+1)$ be the dihedral quandle. Then the third quandle homology $H_3^Q(D_m;\mathbb{Z})$ is $\mathbb{Z}/m\mathbb{Z}$.*

*Proof.* $H_2^Q(D_m;\mathbb{Z}) \cong 0$ is known (see, e.g., [NP2]). Further, $\pi_2(B^Q D_m) \cong \mathbb{Z}/m\mathbb{Z}$ is shown [HN, Corollary 5.9] [3]. Hence, by Theorem 3.9 (i), we have $H_3^Q(D_m;\mathbb{Z}) \cong \mathbb{Z}/m\mathbb{Z}$. □

**Remark 3.12.** Mochizuki [M3] determined the dimension of $H_3^Q(D_m;\mathbb{Z}) \otimes \mathbb{Z}_p$ for any prime $p$. Our result determines the torsion part of $H_3^Q(D_m;\mathbb{Z})$.

---

[3] In [HN], we dealt with a certain group $\Pi_2(D_m)$ instead of $\pi_2(B^Q D_m)$. However, $\pi_2(B^Q D_m) \cong \Pi_2(D_m)$ is known.

# 4 Proof of Theorem 3.1

Our purpose is to prove Theorem 3.1 in §4.2. A computation of $\pi_3(BX)$ will involve a preliminary analysis of the homology $H_i(B_GX; \mathbb{Z})$. For this, in §4.1, we review the rack homology and some topological monoids.

## 4.1 Reviews of rack homology, quandle homology and topological monoid

We first review the rack complexes introduced by [FRS] (see also [EG]) and the quandle homologies defined in [CEGS]. For a quandle $X$, we denote by $C_n^R(X)$ the free $\mathbb{Z}[\mathrm{As}(X)]$-module generated by $n$-elements of $X$. Namely, $C_n^R(X) = \mathbb{Z}[\mathrm{As}(X)]\langle X^n \rangle$. Define a boundary homomorphism $\partial_n : C_n^R(X) \to C_{n-1}^R(X)$, for $n \geq 2$, to be

$$\partial_n(x_1, \ldots, x_n) = \sum_{1 \leq i \leq n} (-1)^i \big( x_i(x_1 * x_i, \ldots, x_{i-1} * x_i, x_{i+1}, \ldots, x_n) - (x_1, \ldots, x_{i-1}, x_{i+1}, \ldots, x_n) \big),$$

and $\partial_1(x_1)$ to be $x_1 - 1 \in \mathbb{Z}[\mathrm{As}(X)]$. Note that the composite $\partial_{n-1} \circ \partial_n$ is zero. Given a left $\mathbb{Z}[\mathrm{As}(X)]$-module $M$, a complex $C_n^R(X; M) = M \otimes_{\mathbb{Z}[\mathrm{As}(X)]} C_n^R(X)$ with $\partial_n$ is called *the rack complex* of $X$ with coefficients $M$. Next, let $C_n^D(X; M)$ be a submodule of $C_n^R(X; M)$ generated by $n$-tuples $(x_1, \ldots, x_n)$ with $x_i = x_{i+1}$ for some $i \in \{1, \ldots, n-1\}$ if $n \geq 2$; otherwise, let $C_1^D(X; M) = 0$. Since $\partial_n(C_n^D(X; M)) \subset C_{n-1}^D(X; M)$, we denote the homology by $H_n^D(X; \mathbb{Z})$. Further, a complex $\big(C_*^Q(X; M), \partial_*\big)$ is defined by the quotient $C_n^R(X; M)/C_n^D(X; M)$. The homology $H_n^Q(X; M)$ is called a *quandle homology* of $X$ with coefficient $M$.

For an $X$-set $Y$, regarding the free module $M = \mathbb{Z}\langle Y \rangle$ as a $\mathbb{Z}[\mathrm{As}(X)]$-module, the complexes $(C_*^R(X; M), \partial_*)$ and $(C_*^Q(X; M), \partial_*)$ are chain isomorphic to the cellular complexes of $B_Y X$ and of $B_Y^Q X$, respectively. As the simplest case, if $Y$ is a single point, then $(C_*^R(X; \mathbb{Z}), \partial_*)$ and $(C_*^Q(X; \mathbb{Z}), \partial_*)$ coincide with the rack complex and the quandle complex of $X$ described in [CJKLS], respectively.

**Remark 4.1.** In the special case $Y = X$, we can identify $C_{n+1}^R(X; \mathbb{Z})$ with $C_n^R(X; \mathbb{Z}\langle X \rangle)$ from the definitions. Moreover, the identification is a chain map, leading $H_{n+1}^R(X; \mathbb{Z}) \cong H_n^R(X; \mathbb{Z}\langle X \rangle) \cong H_n(B_X X; \mathbb{Z})$ (see also [FRS2, Theorem 5.12]).

We next review basic properties of the complexes. Let $X$ be a finite quandle with $\ell$-connected components and $M = \mathbb{Z}$. According to [LN, Theorem 2.2], we have

$$H_1^R(X; \mathbb{Z}) \cong \mathbb{Z}^\ell, \quad H_2^R(X; \mathbb{Z}) \cong H_2^Q(X; \mathbb{Z}) \oplus \mathbb{Z}^\ell, \quad H_3^R(X; \mathbb{Z}) \cong H_3^Q(X; \mathbb{Z}) \oplus H_2^Q(X; \mathbb{Z}) \oplus \mathbb{Z}^{\ell^2}. \quad (5)$$

Further, it is known [CJKS, Lemma 3.2] that

$$\dim_\mathbb{Q}(H_n^R(X; \mathbb{Z}) \otimes \mathbb{Q}) = \ell^n, \qquad \dim_\mathbb{Q}(H_n^Q(X; \mathbb{Z}) \otimes \mathbb{Q}) = \ell \cdot (\ell - 1)^{n-1}. \quad (6)$$

Furthermore, the following long exact sequence splits (see [LN, Theorem 2.1]):

$$\longrightarrow H_{n+1}^Q(X; \mathbb{Z}) \longrightarrow H_n^D(X; \mathbb{Z}) \longrightarrow H_n^R(X; \mathbb{Z}) \longrightarrow H_n^Q(X; \mathbb{Z}) \longrightarrow H_{n-1}^D(X; \mathbb{Z}) \longrightarrow . \quad (7)$$

In particular, we have $H_n^R(X;\mathbb{Z}) \cong H_n^Q(X;\mathbb{Z}) \oplus H_n^D(X;\mathbb{Z})$.

Given a connected quandle $X$ and an $X$-set $Y$, we observe the rack space $B_Y X$ as a covering as follows. Note that the classing map $Y \to \{\text{pt.}\}$ induces a continuous map $B_Y X \to BX$. It is known (see [FRS, Theorem 3.7]) that this is a covering; hence, so is the induced map $B_Y^Q X \to B^Q X$. As a special case, if $Y = G$ is a quotient group of $\mathrm{As}(X)$ (hence $B_G X$ is path-connected), then the two coverings $B_G X \to BX$ and $B_G^Q X \to B^Q X$ are principal $G$-bundles over $BX$ and $B^Q X$, respectively (see [Cla, Proposition 6]). Recalling $\pi_1(BX) \cong \mathrm{As}(X)$ from §2.1, the group $\pi_1(B_G X)$ is exactly the kernel of $p_G : \mathrm{As}(X) \to G$

The action of $\pi_1(B_G X)$ on the higher homotopy groups $\pi_*(B_G X)$ is known to be trivial (see [FRS2, Proposition 5.2] and [Cla, Proposition 10]). Furthermore, the action of $\pi_1(BX)$ on the higher homology group $H_*(B_G X)$ is also trivial (see [EG, Lemma 3.1] or [Cla, Remark 7]); hence, the action on $H_*(B_G^Q X)$ is also trivial, since $\pi_1(BX) = \pi_1(B^Q X)$ and $H_*(B_G^Q X)$ is a quotient of $H_*(B_G X)$ by (7).

Finally, we discuss topological monoids on some rack spaces introduced by Clauwens [Cla]. We consider a natural map $X \to \mathrm{Inn}(X)$ given by $x \mapsto (\bullet * x)$, which yields an epimorphism $\mathrm{As}(X) \to \mathrm{Inn}(X)$. Let $G$ be a quotient group of $\mathrm{As}(X)$. We here the epimorphism induces $G \to \mathrm{Inn}(X)$. Then, Clauwens introduced an operation given by

$$\mu : (G \times [0,1]^n \times X^n) \times (G \times [0,1]^m \times X^m) \to G \times [0,1]^{n+m} \times X^{n+m},$$

$$\mu([g; t_1, \ldots, t_n, x_1, \ldots, x_n], [h; t_1', \ldots, t_m', x_1', \ldots, x_m']) :=$$

$$[gh; t_1, \ldots, t_n, t_1', \ldots, t_m', x_1 \cdot h, \ldots, x_n \cdot h, x_1', \ldots, x_m'],$$

where $n, m \in \mathbb{Z}_{\geq 0}$; this further gives rise to a topological monoid structure on $B_G X$ (see [Cla, §2.5]). In particular, $\pi_1(B_G X)$ is an abelian group. Since $B_G X$ is a so-called nilpotent space, we often deal with the localization of $B_G X$ at a prime $p$.

## 4.2 Proof of Theorem 3.1

*Proof.* Consider $G = \mathrm{Inn}(X)$ which is a finite group. Since the rack space $B_G X$ is a topological monoid and is a path-connected CW complex of finite type, $B_G X$ is homotopic to a (based) loop space of some simply connected CW complex $W$ of finite type, that is, $B_G X \simeq \Omega W$ as an $H$-space (see [Mil, Theorem 1.5 and §2] for details). In particular, it is without saying that $\pi_*(B_G X)$ is finitely generated; hence, so is the quotient $\pi_3^Q(BX)$.

Next, we discuss rational homologies of $B_G X$ and $B_G^Q X$. Recall that the projections $p_G^Q : B_G^Q X \to B^Q X$ and $p_G : B_G X \to BX$ are coverings of degree $|G|$, and that the actions of $\pi_1(BX)$ on $H_*(B_G^Q X)$ and on $H_*(B_G X)$ are trivial. It then follows from the transfer maps that $p_G^Q$ and $p_G$ induce two isomorphisms:

$$(p_G)_* : H_n(B_G X; \mathbb{Q}) \cong H_n(BX; \mathbb{Q}), \quad (p_G^Q)_* : H_n(B_G^Q X; \mathbb{Q}) \cong H_n(B^Q X; \mathbb{Q}). \tag{8}$$

Next, we discuss $\pi_*(B_GX) \otimes \mathbb{Q}$. Since $B_GX \simeq \Omega W$, we recall the known following formula (see [FHT, §33(c)]):

$$\sum_{i \geq 0} \dim(H^i(B_GX; \mathbb{Q})) \, t^i = \prod_{i \geq 0} \frac{(1+t^{2i+1})^{r_{2i+1}}}{(1-t^{2i})^{r_{2i}}} \in \mathbb{Z}[[t]],$$

where we put $r_i = \dim(\pi_i(B_GX) \otimes \mathbb{Q})$. By (6) and (8), we then easily have $r_2 = (\ell^2+\ell)/2$ and $r_3 = (\ell^3-\ell)/3$. Since $\pi_2(BX) \cong \pi_2(B^QX) \oplus \mathbb{Z}^\ell$ is known (see [N1, Proposition 3.12]), we conclude $\dim(\pi_2(B^QX) \otimes \mathbb{Q}) = (\ell^2-\ell)/2$.

We next calculate $\pi_3^Q(BX) \otimes \mathbb{Q}$. For this, we equip $H^*(B_GX; \mathbb{Q})$ with a Hopf algebra arising from the monoid structure, and denote by $P^*(B_GX)$ the set of primitive elements of $H^*(B_GX; \mathbb{Q})$. By Milnor-Moore theorem (see, e.g., [FHT, Theorem 21.5]), $\pi_*(B_GX) \otimes \mathbb{Q} \cong P^*(B_GX)$ and the dual Hurewicz homomorphism $(\mathfrak{H})^* \otimes \mathbb{Q} : H^*(B_GX; \mathbb{Q}) \to \mathrm{Hom}(\pi_*(B_GX), \mathbb{Q})$ coincides with the projection $H^*(B_GX; \mathbb{Q}) \to P^*(B_GX)$.

We now explain the isomorphism (10) below. For this, by (8) and the splitting in (7), we notice that the induced map $(i_X)_* : H_*(B_GX; \mathbb{Q}) \to H_*(B_G^QX; \mathbb{Q})$ is a splitting surjection, where $i_X$ is the inclusion $B_GX \hookrightarrow B_G^QX$. Further, we consider the long exact sequences of cohomology and homotopy groups of the pair $(B_G^QX, B_GX)$, and have the following commutative diagram:

$$\begin{array}{ccccccc}
0 & \longrightarrow & H^3(B_G^QX; \mathbb{Q}) & \xrightarrow{(i_X)^*_{H^3}} & H^3(B_GX; \mathbb{Q}) & \xrightarrow{\partial^*} & H^4(B_G^QX, B_GX; \mathbb{Q}) \\
 & & \mathfrak{H}_3^Q \downarrow & & \mathfrak{H}_3^R \downarrow & & \\
0 & \longrightarrow & \mathrm{Hom}(\pi_3^Q(BX), \mathbb{Q}) & \xrightarrow{(i_X)^*} & \mathrm{Hom}(\pi_3(B_GX), \mathbb{Q}) & &
\end{array} \quad (9)$$

Here the top and bottom sequences are exact, and the vertical arrows are the Hurewicz homomorphisms. By Milnor-Moore theorem again, we obtain

$$\pi_3^Q(BX) \otimes \mathbb{Q} \cong P^3(B_GX) \cap \mathrm{Ker}\big(\partial^* : H^3(B_GX; \mathbb{Q}) \to H^4(B_G^QX, B_GX; \mathbb{Q})\big). \quad (10)$$

we will calculate $P^3(B_GX) \cap \mathrm{Ker}(\partial^*)$ using the complexes in §4.1. Inspired by the monoid structure on $B_GX$, we consider a map $\mu : (G \times X^n) \times (G \times X^m) \to G \times X^{n+m}$ given by

$$\mu\big((g; x_1, \ldots, x_n), (h; x'_1, \ldots, x'_m)\big) := (gh; x_1 \cdot h, \ldots, x_n \cdot h, x'_1, \ldots, x'_m).$$

This provides $H_D^*(X; \mathbb{Q}[G])$ and $H_R^*(X; \mathbb{Q}[G])$ with Hopf algebra structures[4]. Further, note that the inclusion $i_G^D : C_*^D(X; \mathbb{Q}[G]) \to C_*^R(X; \mathbb{Q}[G])$ induces a Hopf algebra homomorphism $(i_G^D)^* : H_R^*(X; \mathbb{Q}[G]) \to H_D^*(X; \mathbb{Q}[G])$. Let us identify the above map $\partial^*$ with the induced map $(i_G^D)^*$ by definitions. By (7) we then notice that

$$\mathrm{Ker}(\partial^*) \cap H^3(B_GX; \mathbb{Q}) \cong \mathrm{Ker}((i_G^D)^*) \cap H_R^3(X; \mathbb{Q}[G]) \cong H_Q^3(X; \mathbb{Q}) \cong \mathbb{Q}^{\ell(\ell-1)^2},$$

$$\mathrm{Ker}((i_G^D)^*) \cap P^1(B_GX) \cong H_R^1(X; \mathbb{Q}[G]) \cong \mathbb{Q}^\ell.$$

---
[4]Here the products are defined by the cup products. See [Cla, §2.8] for the explicit formula of the cup product.

In conclusion, the dimension of $P^3(B_G X) \cap \mathrm{Ker}(\partial^*)$ in (10) is equal to

$$\dim\bigl(H^3_Q(X;\mathbb{Q}[G])\bigr) - \dim\bigl(\Lambda^3\bigl(P^1(B_G X)\bigr)\bigr) - \dim\bigl(P^1(B_G X) \otimes P^2(B_G X)\bigr)$$

$$= \ell(\ell-1)^2 - \frac{\ell(\ell-1)(\ell-2)}{6} - l\frac{(\ell^2+\ell)}{2} = \frac{\ell(\ell-1)(\ell-2)}{3}. \qquad \square$$

**Remark 4.2.** We roughly explain a topological meaning of the quandle homotopy invariant after tensoring with $\mathbb{Q}$ as follows. To illustrate, we put a collapsing map $X \to T_\ell$ on each connected components, where $T_\ell$ is the trivial quandle of order $\ell$. We see that this map induces $\pi_3^Q(BX) \otimes \mathbb{Q} \cong \pi_3^Q(BT_\ell) \otimes \mathbb{Q}$ by the discussion in the proof; we may assume $X = T_\ell$. Then we can obtain all cocycles of $H^3_R(T_\ell;\mathbb{Q}) \cong \mathbb{Q}^{l^3}$, since the coboundary maps $\delta_n$ are zero by definition. Since $\pi_3^Q(BT_\ell) \otimes \mathbb{Q}$ is derived from the primitive elements of $H^R_3(T_\ell;\mathbb{Q})$ by Milnor-Moore theorem, we can evaluate elements of $\pi_3^Q(BT_\ell) \otimes \mathbb{Q}$ by the cocycles.

On the other hand, we consider the bordism group $L^2_\ell$ with $\ell$-connected components (see [CKSS, §1] for the definition). An isomorphism $L^2_\ell \otimes \mathbb{Q} \cong \mathbb{Q}^{\frac{\ell(\ell-1)(\ell-2)}{3}}$ is known; further it is shown [CKSS] that $L^2_\ell \otimes \mathbb{Q}$ is generated by "Hopf 2-links". Recall the map $\Xi_X(D; \bullet) : \mathrm{Col}_X(D) \to \pi_3(B^Q X)$ with $X = T_\ell$ in §2.2. By running over all $T_\ell$-colorings of all broken diagrams, the maps give rise to a homomorphism $L^2_\ell \otimes \mathbb{Q} \to \pi_3^Q(BX) \otimes \mathbb{Q}$. Since we easily compute the rational quandle homotopy invariant of the Hopf 2-links with $X = T_\ell$ by pairing with the previous cocycles, the homomorphism turns to be an isomorphism: $L^2_\ell \otimes \mathbb{Q} \cong \pi_3^Q(BX) \otimes \mathbb{Q}$.

## 5 Proofs of Theorems 3.5 and 3.9.

Our purpose in this section is to prove Theorem 3.5 in §5.1, and Theorem 3.9 in §5.2.

We now outline the proof of Theorem 3.5 as follows. To see this, we first aim to study $\pi_3(BX)$, since $\pi_3^Q(BX)$ is a quotient of $\pi_3(BX)$. Let $\widetilde{BX}$ be the universal covering of $BX$. By the monoid structure of $BX$ explained in §4.1, $\widetilde{BX}$ is a loop space of a 2-connected CW-complex. Thanks to the fact [AP] that the first $k$-invariant of $\widetilde{BX}$ is annihilated by 2, we thus have $\pi_3(BX) = \pi_3(\widetilde{BX}) \cong H_3(\widetilde{BX};\mathbb{Z}) \oplus \pi_3(\Omega S^2)$ by a standard argument of the Postnikov tower. In §5.1, we prove $H_3(\widetilde{BX};\mathbb{Z}) \cong H^Q_4(X;\mathbb{Z})$ using techniques of quandle homologies (see (11)). Finally, we show that $H^Q_4(X;\mathbb{Z})$ does not vanish in $\pi_3^Q(BX)$, leading $\pi_3^Q(BX) \cong H^Q_4(X;\mathbb{Z})$ as desired. Theorem 3.9 is also proven in a similar way.

To carry the outline, we here prepare two propositions, which is proven in §5.3.

**Proposition 5.1.** *Let $X$ be a regular Alexander quandle of finite order. Let $G = \mathrm{Inn}(X)$. Then $H_2(B_G X;\mathbb{Z}) \cong H^Q_3(X;\mathbb{Z}) \oplus H^Q_2(X;\mathbb{Z}) \oplus \mathbb{Z}$. Furthermore, if $H^Q_2(X;\mathbb{Z}) \cong 0$, then*

(i) $H_3(B_G X;\mathbb{Z}) \cong H^Q_4(X;\mathbb{Z}) \oplus H^Q_3(X;\mathbb{Z}) \oplus \mathbb{Z}$.

(ii) $H_3(B^Q_G X;\mathbb{Z}) \cong H^Q_4(X;\mathbb{Z}) \oplus H^Q_3(X;\mathbb{Z})$.

(iii) The map $H_3(B_G X; \mathbb{Z}) \to H_3(B_G^Q X; \mathbb{Z})$ induced by the inclusion $B_G X \hookrightarrow B_G^Q X$ coincides with the projection in the sense of the above presentations (i) and (ii).

**Proposition 5.2.** *Let $X$ be a regular Alexander quandle of finite order. Let $\mathrm{As}(X) \to \mathrm{Inn}(X)$ be the epimorphism in §4.1. Then the kernel is isomorphic to $\pi_1(B_G X) = H_2^R(X; \mathbb{Z}) = \mathbb{Z} \oplus H_2^Q(X; \mathbb{Z})$.*

Furthermore, we require an elementary lemma:

**Lemma 5.3.** *Let $\mathcal{M}$ be a topological monoid with the binary operation $\mu$. If $\mathcal{M}$ is a path-connected CW-complex and $\pi_1(\mathcal{M}) \cong \mathbb{Z}$, then $\mathcal{M}$ is homotopic to $S^1 \times \widetilde{\mathcal{M}}$. Here $\widetilde{\mathcal{M}}$ is the universal covering of $\mathcal{M}$.*

*Proof.* Choose a representative $f : S^1 \to \mathcal{M}$ of a generator of $\pi_1(\mathcal{M})$. Let $g : \widetilde{\mathcal{M}} \to \mathcal{M}$ denote the universal covering map. The composite map $\mu \circ (f \times g) : S^1 \times \widetilde{\mathcal{M}} \to \mathcal{M}$ turns out to give a weak homotopy equivalence as desired. □

## 5.1 Proof of Theorem 3.5

*Proof.* From now on, $X$ is assumed to be a regular Alexander quandle with $H_2^Q(X; \mathbb{Z}) \cong 0$. By Proposition 5.2, the kernel of the epimorphism $\mathrm{As}(X) \to \mathrm{Inn}(X)$ is entirely $\mathbb{Z}$. Hence, putting $G = \mathrm{Inn}(X)$, we have $B_G X \simeq S^1 \times \widetilde{BX}$ by Lemma 5.3; thus it follows from Proposition 5.1 and Kunneth formula that

$$H_2(\widetilde{BX}; \mathbb{Z}) \simeq \mathbb{Z} \oplus H_3^Q(X; \mathbb{Z}), \qquad H_3(\widetilde{BX}; \mathbb{Z}) \simeq H_4^Q(X; \mathbb{Z}). \tag{11}$$

Let us calculate the localizations of $\pi_3(\widetilde{BX})$ at every prime $p$. Since $\widetilde{BX}$ is also a topological monoid; $\widetilde{BX}$ is a loop space. Let $\widetilde{BX}_i$ denote the $i$-th stage of Postnikov tower. It immediately follows from [AP, Theorem 3.2] that the first $k$-invariant $k^4 \in H^4(K(\pi_2(\widetilde{BX}), 2); \pi_3(\widetilde{BX}))$ satisfies $2k^4 = 0$ (see also [Sou, Proposition 3]). Hence, since the order $|X|$ is odd, the localization of $\widetilde{BX}_3$ at a prime $p > 2$ is reduced to be

$$(\widetilde{BX}_{(p)})_3 \simeq K(\pi_2(BX)_{(p)}, 2) \times K(\pi_3(BX)_{(p)}, 3). \tag{12}$$

Hence, the Hurewicz map $\mathfrak{H}_{(p)}$ passes to an isomorphism $\pi_3(BX)_{(p)} \cong H_3(\widetilde{BX}; \mathbb{Z})_{(p)}$, since $H_3(K(\pi_2(BX), 2); \mathbb{Z})$ vanishes (see [McC, §$8^{bis}$]). As a consequence, combining with (11) immediately gives $\pi_3(BX)_{(p)} \cong H_3(\widetilde{BX}; \mathbb{Z})_{(p)} \cong H_4^Q(X; \mathbb{Z})_{(p)}$.

We consider the last case $p = 2$. It is known [FRS2, Theorem 5.12] that $BT_1$ is homotopic to $\Omega S^2$, where $T_1$ is the trivial quandle of order 1. Let $X \to T_1$ be the collapsing map. This induces a topological monoid homomorphism $f : B_G X \to \Omega S^2$. By Lemma 5.7 below, the torsion part of $H_*(B_G X; \mathbb{Z})$ is annihilated by $|X|$; thereby, $(f_{(2)})_* : H_*(B_G X; \mathbb{Z})_{(2)} \to H_*(\Omega S^2; \mathbb{Z})_{(2)} \ (\cong \mathbb{Z}_{(2)}[\zeta])$ is isomorphic, where $\zeta$ is a generator of $H_1(\Omega S^2; \mathbb{Z}) \cong \mathbb{Z}$. Hence, by the Whitehead theorem, the cellular map $f$ gives rise to a homotopy equivalence $f_{(2)} : B_G X_{(2)} \to \Omega S^2_{(2)}$. In particular, $\pi_3(B_G X_{(2)}) \cong \pi_3(\Omega S^2)_{(2)} \cong \pi_4(S^2)_{(2)} \cong \mathbb{Z}/2\mathbb{Z}$. Further, by Theorem 3.1 and (6), we note $\pi_3(BX)_{(0)} \cong H_4^Q(X; \mathbb{Z})_{(0)} \cong 0$. In summary, we conclude $\pi_3(BX) \cong H_4^Q(X; \mathbb{Z}) \oplus \mathbb{Z}/2\mathbb{Z}$ as required.

Finally, we discuss $\pi_3^Q(BX)$. Putting the covering $g : \widetilde{BX} \to B_G X$ and the inclusion $i_X : B_G X \hookrightarrow B_G^Q X$, we have a commutative diagram.

$$\begin{array}{ccccc}
\pi_3(\widetilde{BX})_{(p)} & \xrightarrow{\sim} & \pi_3(B_G X)_{(p)} & \xrightarrow{(i_X)_*} & \pi_3^Q(BX)_{(p)} \\
\downarrow \widetilde{\mathfrak{H}}_3 & & \downarrow \mathfrak{H}_3^R & & \downarrow \mathfrak{H}_3^Q \\
H_3(\widetilde{BX}; \mathbb{Z})_{(p)} & \xrightarrow{g_*} & H_3(B_G X; \mathbb{Z})_{(p)} & \xrightarrow{(i_X)_*} & H_3(B_G^Q X; \mathbb{Z})_{(p)}
\end{array}$$

For $p > 2$, we remark that the localization of the Hurewicz map $\widetilde{\mathfrak{H}}_3$ is isomorphic by (12). By Lemma 5.3, the map $g_*$ is injective; hence, so is $\mathfrak{H}_3^R$. We note that, by Proposition 5.1, the above map $(i_X)_* : H_3(B_G X; \mathbb{Z}) \to H_3(B_G^Q X; \mathbb{Z})$ is isomorphic. Thereby the bottom map $(i_X)_* : \pi_3(B_G X)_{(p)} \to \pi_3(B_G^Q X)_{(p)}$ is an isomorphism.

To complete the proof, it suffices to show that the direct summand $\mathbb{Z}/2\mathbb{Z}$ of $\pi_3(BX)$ is sent to the zero via $(i_X)_* : \pi_3(BX) \to \pi_3(B^Q X)$. By the above discussion of the direct summand $\mathbb{Z}/2\mathbb{Z}$, we may assume $X = T_1$. It is known [San, CKSS] that $\pi_3(BT_1) \cong \mathbb{Z}/2\mathbb{Z}$ is isomorphic to the framed link bordism group ($\cong \mathbb{Z}/2\mathbb{Z}$), and that the generator of the bordism group is represented by an embedding of a pair of unknotted tori (see [San, Example 1.12] for details). Hence, from the definition of $B^Q X$, for any $T_1$-coloring $C$ of the tori, the cellular map $\xi_{D,C}$ in §2.2 is bounded by a 3-cell of the subspace $X_Y^D$ in §2.1. Namely, the homotopy class $[\xi_{D,C}]$ vanishes in $\pi_3(B^Q X)$, which completes the proof. □

**Remark 5.4.** In a similar manner, we can calculate the forth homotopy groups $\pi_4(BX)$ of regular Alexander quandles $X$ of odd order with $H_2^Q(X; \mathbb{Z}) = 0$ as follows. In fact, it is shown [AP, Theorem 4.6] that the second $k$-invariant of $\widetilde{BX}$ is annihilated by 16; hence the forth stage $(\widetilde{BX})_4$ localized at $p > 2$ is presented by

$$(\widetilde{BX}_{(p)})_4 \simeq K(\pi_2(BX)_{(p)}, 2) \times K(\pi_3(BX)_{(p)}, 3) \times K(\pi_4(BX)_{(p)}, 4). \tag{13}$$

About the case $p = 2$, it follows from the proceeding map $f_{(2)}$ that $\pi_4(BX)_{(2)} = \pi_4(\Omega S^2)_{(2)} = \pi_5(S^2)_{(2)} = \mathbb{Z}_2$. Thereby, noting $H_4((\widetilde{BX})_4; \mathbb{Z}) \cong H_4(\widetilde{BX}; \mathbb{Z})$, Lemma 5.3 includes that the group $\pi_4(BX)$ can be calculated from $\pi_2(BX)$, $\pi_3(BX)$ and $H_4(B_G X; \mathbb{Z})$.

For instance, let us compute $\pi_4(BX)$ of $X = \mathbb{Z}[T]/(p, T+1)$ with $p > 2$. Recall that $\pi_3(BX) = \mathbb{Z}_{2p}$ by Corollary 3.7, and that $\pi_2(BX) = \mathbb{Z} \oplus \mathbb{Z}_p$ by Theorem 3.9. In addition, $H_4(B_G X; \mathbb{Z}) = \mathbb{Z} \oplus \mathbb{Z}_p^3$ is known [Cla, §1.4], which means $H_4(\widetilde{BX}; \mathbb{Z}) = \mathbb{Z} \oplus \mathbb{Z}_p^2$ by Lemma 5.3. Therefore, by (13) we have $\pi_4(BX) = \mathbb{Z}_2$ (cf. Corollary 3.8).

## 5.2 Proofs of Theorem 3.9 and Proposition 3.4.

We let $X$ be a regular Alexander quandle without the assumption $H_2^Q(X; \mathbb{Z}) = 0$. We fix $G = \text{Inn}(X)$. To prove Theorem 3.9, we use the exact sequence in [McC, Lemma $8^{bis}.27$] in the case of $n = 2$: Namely, for the rack space $BX$, the exact sequence is written in

$$H_3(\pi_1(B_GX); \mathbb{Z}) \longrightarrow \pi_2(B_GX) \xrightarrow{\mathfrak{H}} H_2(B_GX; \mathbb{Z}) \longrightarrow H_2(\pi_1(B_GX); \mathbb{Z}) \longrightarrow 0. \qquad (14)$$

Here $H_*(\pi_1(B_GX); \mathbb{Z})$ is the group homology of $\pi_1(B_GX)$.

*Proof of Theorem 3.9.* (i) By assumption of $H_2^Q(X; \mathbb{Z}) = 0$, we have $\pi_1(B_GX) \cong \mathbb{Z}$ by Proposition 5.2. Noting $H_r(\mathbb{Z}; \mathbb{Z}) \cong 0$ for $r \geq 2$, the Hurewicz homomorphism $\mathfrak{H}$ in (14) is isomorphic. Recall $H_2(B_GX; \mathbb{Z}) \cong H_3^Q(X; \mathbb{Z}) \oplus H_2^Q(X; \mathbb{Z}) \oplus \mathbb{Z}$ by Proposition 5.1. Since $\pi_2(BX) \cong \mathbb{Z} \oplus \pi_2(B^QX)$, we have $\pi_2(BX) \cong H_3^Q(X; \mathbb{Z})$ as required.

(ii) To obtain the required sequence (4), we will show that the map $\mathfrak{H}$ in (14) is a splitting injection. Recall that $B_GX$ is a path-connected loop space. It immediately follows from [AP, Theorem 3.2] that the first $k$-invariant $k^3 \in H^3(\pi_1(B_GX); \pi_2(B_GX))$ satisfies $2k^3 = 0$. By Proposition 5.2, $\pi_1(B_GX) \cong H_2^R(X; \mathbb{Z})$, which implies that the group cohomology $H^3(\pi_1(B_GX); \pi_2(B_GX))$ is annihilated by $|X|$. Hence, the first $k$-invariant is zero. That is, the second stage of the Postnikov tower is homotopic to $K(\pi_1(B_GX); 1) \times K(\pi_2(B_GX); 2)$, and the map $\mathfrak{H}$ is a splitting injection (cf. [AP, Page 3]).

We now calculate each terms in the sequence (14). Using Proposition 5.2, the second group homology in (14) is expressed as

$$H_2(\pi_1(B_GX); \mathbb{Z}) \cong H_2(\mathbb{Z} \oplus H_2^Q(X; \mathbb{Z}); \mathbb{Z}) \cong H_2^Q(X; \mathbb{Z}) \oplus \left( H_2^Q(X; \mathbb{Z}) \wedge_{\mathbb{Z}} H_2^Q(X; \mathbb{Z}) \right).$$

Since $H_2(B_GX; \mathbb{Z}) = \mathbb{Z} \oplus H_2^Q(X; \mathbb{Z}) \oplus H_3^Q(X; \mathbb{Z})$ by Proposition 5.1, the sequence (14) is rewritten,

$$0 \to \pi_2(BX) \xrightarrow{\mathfrak{H}} \mathbb{Z} \oplus H_2^Q(X; \mathbb{Z}) \oplus H_3^Q(X; \mathbb{Z}) \longrightarrow H_2^Q(X; \mathbb{Z}) \oplus \left( H_2^Q(X; \mathbb{Z}) \wedge_{\mathbb{Z}} H_2^Q(X; \mathbb{Z}) \right) \to 0.$$

Since $\pi_2(BX) \cong \mathbb{Z} \oplus \pi_2(B^QX)$ and the restriction of $\mathfrak{H}$ to the $\mathbb{Z}$-part is isomorphic (see [N1, Proposition 3.12]), we immediately obtain the required sequence (4). □

**Remark 5.5.** As is seen in the proof, the Hurewicz homomorphism $\mathfrak{H}$ in the sequence (14) is injective. Hence, via $\mathfrak{H}$, any element of $\pi_2(B^QX)$ can be evaluated by some cocycles $\psi \in H_Q^2(X; A) \oplus H_Q^3(X; A)$ with trivial coefficients. In conclusion, there exist $\psi_1, \ldots, \psi_m \in H_Q^2(X; A) \oplus H_Q^3(X; A)$ such that any quandle cocycle invariant is a linear sum of the quandle cocycle invariants associated with $\psi_1, \ldots, \psi_m$.

Furthermore, by a similar discussion, we show Proposition 3.4 as follows:

*Proof of Proposition 3.4.* Let $Y_i$ denote the $i$-th stage of the Postnikov tower of $B_GX$. We recall the following exact sequence in [McC, Lemma $8^{bis}.27$].

$$H_4(Y_2; \mathbb{Z}) \longrightarrow \pi_3(B_GX) \xrightarrow{\mathfrak{H}} H_3(B_GX; \mathbb{Z}) \longrightarrow H_3(Y_2; \mathbb{Z}) \longrightarrow 0. \qquad (15)$$

Let us calculate the each terms in (15). By Lemma 5.9, $H_3(B_GX; \mathbb{Z})$ is isomorphic to the rack homology $H_4^R(X; \mathbb{Z})$. Moreover, it follow from the previous proof of Theorem 3.9 that $Y_2$ is homotopic to $K(\pi_1(B_GX), 1) \times K(\pi_2(B_GX), 2)$. Further, recall the isomorphism $H_2^R(X; \mathbb{Z}) \cong \pi_1(B_GX)$ from Proposition 5.2. Therefore, we easily obtain the required sequence from the sequence (15) and the Kunneth formula. □

## 5.3 Proofs of Propositions 5.1 and 5.2

In §5.3 and 5.4, $X$ is assumed to be a finite connected Alexander quandle. We fix $G = \mathrm{Inn}(X)$. Recall $H_n(B_G X; \mathbb{Z}) \cong H_n^R(X; \mathbb{Z}[\mathrm{Inn}(X)])$ in §4.1. To show Proposition 5.1, we first study $\mathrm{Inn}(X)$.

**Lemma 5.6.** *Let $e \in \mathbb{Z}_{\geq 0}$ be the minimal number satisfying $(T^e - 1)X = 0$. Then $\mathrm{Inn}(X) \cong \mathbb{Z}/e\mathbb{Z} \ltimes X$. Here $\mathbb{Z}/e\mathbb{Z}$ acts on $X$ by the multiplication by $T$.*

*Proof.* For $(\varepsilon, x) \in \mathbb{Z}/e\mathbb{Z} \ltimes X$, we put an element of $\mathrm{Map}(X, X)$ sending $y \in X$ to $T^\varepsilon y + (1 - T)x \in X$. This then gives rise to an epimorphism $\mathbb{Z}/e\mathbb{Z} \ltimes X \to \mathrm{Inn}(X)\ (\subset \mathfrak{S}_{|X|})$. We claim that this is injective. Actually, if $(\varepsilon, x) \in \mathbb{Z}/e\mathbb{Z} \ltimes X$ satisfies $y = T^\varepsilon y + (1 - T)x$ for any $y \in X$, then we easily obtain $(\varepsilon, x) = (0, 0)$. □

As a result, we see that the surjection $\mathrm{As}(X) \to \mathrm{Inn}(X)$ in §4.1 coincides with a homomorphism $\mathrm{As}(X) \to \mathbb{Z}/e\mathbb{Z} \ltimes X$ sending $x$ to $(1, x)$. We also observe another lemma:

**Lemma 5.7.** *Let $X$ be a finite connected Alexander quandle. Let $M$ be $\mathbb{Z}[\mathrm{Inn}(X)]$. Then the torsion subgroup of $H_*^R(X; M) = H_*(B_G X; \mathbb{Z})$ is annihilated by $|X|$.*

We defer its proof until §5.4. We next study a relation between $H_n^Q(X; M)$ and $H_n(B_G^Q X; \mathbb{Z})$. To see this, we review a complex $C_n^L(X)$ introduced in [LN, §2]. This $C_n^L(X)$ is defined to be a subcomplex of $C_n^D(X)$ generated by $n$-tuples $(x_1, \ldots, x_n) \in X^n$ with $x_i = x_{i+1}$ for some $2 \leq i < n$. It is shown [LN, Lemma 9] that there is an isomorphism of chain complexes:

$$C_*^D(X; \mathbb{Z}) \cong C_{*-1}^Q(X; \mathbb{Z}) \oplus C_*^L(X; \mathbb{Z}), \tag{16}$$

where the direct summand $C_*^L(X; \mathbb{Z})$ is obtained from the inclusion $C_*^L(X; \mathbb{Z}) \hookrightarrow C_*^D(X; \mathbb{Z})$.

**Lemma 5.8.** *If $H_2^Q(X; \mathbb{Z}) = 0$, then $H_4^L(X; \mathbb{Z}) \cong \mathbb{Z}$.*

We later prove this in §5.4. Meanwhile, we consider a complex $C_*^R(X; \mathbb{Z}\langle X \rangle)$ in Remark 4.1. Under the identification $C_n^R(X; \mathbb{Z}\langle X \rangle) \cong C_{n+1}^R(X; \mathbb{Z})$, as the restrictions, we have a chain isomorphism $C_n^D(X; \mathbb{Z}\langle X \rangle) \cong C_{n+1}^L(X; \mathbb{Z})$. In conclusion, by (7) and (16), we thus have

$$H_n(B_X^Q X; \mathbb{Z}) \cong H_n^Q(X; \mathbb{Z}\langle X \rangle) \cong H_{n+1}^Q(X; \mathbb{Z}) \oplus H_n^Q(X; \mathbb{Z}). \tag{17}$$

**Lemma 5.9.** *Let $X$ be a regular Alexander quandle of finite order. Then $H_n(B_G X; \mathbb{Z}) \cong H_{n+1}^R(X; \mathbb{Z})$ for $n \geq 1$. Furthermore, $H_n(B_G^Q X; \mathbb{Z}) \cong H_{n+1}^Q(X; \mathbb{Z}) \oplus H_n^Q(X; \mathbb{Z})$.*

*Proof.* Let us regard $X$ as an $X$-set (see §2.1). We consider the two maps $\mathrm{Inn}(X) \to X$ and $X \to \{\mathrm{pt.}\}$. They then give rise to two coverings $B_G X \to B_X X$ and $B_X X \to BX$ (cf. [Cla, §4.1]). Further, note that the covering $B_G X \to B_X X$ is of degree $e$ by Lemma 5.6. Hence, the transfer map yields an isomorphism $H_n(B_G X; \mathbb{Z})_{(p)} \cong H_n(B_X X; \mathbb{Z})_{(p)}$ where the prime $p$ divides $|X|$ (cf. [Cla, Propositions 22 and 23]). Furthermore, it is evident that $H_n(B_G X; \mathbb{Z})_{(0)} \cong H_n(BX; \mathbb{Z})_{(0)} \cong \mathbb{Q}$ by (6). It is shown [N2, Theorem 6.3] that the torsion subgroup of $H_n(BX;\mathbb{Z})$ is annihilated by $|X|$. Combining with Lemma 5.7, we have $H_n(B_GX;\mathbb{Z}) \cong H_n(B_XX;\mathbb{Z}) \cong H_{n+1}(BX;\mathbb{Z}) \cong H^R_{n+1}(X;\mathbb{Z})$.

To prove the latter part, similarly, we can see $H_*(B^Q_GX;\mathbb{Z}) \cong H_*(B^Q_XX;\mathbb{Z})$. Therefore, we have $H_n(B^Q_GX;\mathbb{Z}) \cong H^Q_{n+1}(X;\mathbb{Z}) \oplus H^Q_n(X;\mathbb{Z})$ by (17), which completes the proof. $\square$

We now prove Propositions 5.1 and 5.2.

*Proof of Proposition 5.1.* The desired isomorphism $H_2(B_GX;\mathbb{Z}) \cong H^Q_3(X;\mathbb{Z}) \oplus H^Q_2(X;\mathbb{Z}) \oplus \mathbb{Z}$ immediately follows from Lemma 5.9 and (7). Similarly, (ii) is follows from the lemma.

To prove (i), using (7) and (16), Lemma 5.9 results in

$$H_3(B_GX;\mathbb{Z}) \cong H^R_4(X;\mathbb{Z}) \cong H^Q_4(X;\mathbb{Z}) \oplus H^Q_3(X;\mathbb{Z}) \oplus H^L_4(X;\mathbb{Z}). \tag{18}$$

Therefore, the required isomorphism (i) is obtained from Lemma 5.8. Finally, the proof of (iii) immediately follows from the sequence (7) and the isomorphism (17). $\square$

*Proof of Proposition 5.2.* Recall that the kernel of the epimorphism $\mathrm{As}(X) \to \mathrm{Inn}(X) = G$ is $\pi_1(B_GX)$ and abelian (see §4.1). Hence, $\pi_1(B_GX) \cong H_1(B_GX;\mathbb{Z}) \cong H^R_2(X;\mathbb{Z})$ by Lemma 5.9. $\square$

## 5.4 Proofs of Lemmas 5.8 and 5.7

Before proving Lemma 5.8, we remark the condition of $H^Q_2(X;\mathbb{Z}) \cong 0$. By (5), we have $H^R_2(X;\mathbb{Z}) \cong \mathbb{Z}$. It is known [CJKS, LN] that, for any $x \in X$, the 2-cycle of the form $(x,x) \in C^R_2(X;\mathbb{Z})$ is a generator of $H^R_2(X;\mathbb{Z}) \cong \mathbb{Z}$. Namely, any cycle of $C^R_2(X;\mathbb{Z})$ is homologous to $\alpha(x,x)$ for some $\alpha \in \mathbb{Z}$.

*Proof of Lemma 5.8.* Consider a cycle $\sigma \in C^L_4(X;\mathbb{Z})$, i.e., $\partial_4(\sigma) = 0$. It is enough to show that the cycle $\sigma$ is homologous to a sum of $(y,y,y,y)$s for some $y \in X$. Let us expand $\sigma$ as $\sum_i \alpha_i(a_i,b_i,b_i,c_i) + \sum_j \beta_j(d_j,e_j,f_j,f_j)$ for some $a_i,b_i,c_i,d_j,e_j,f_j \in X$ and $\alpha_i,\beta_j \in \mathbb{Z}$, where $e_j \neq f_j$.

We show (19) below. For any $g_i \in X$, we note that

$$\partial_5(a_i,g_i,b_i,b_i,c_i) = (a_i,b_i,b_i,c_i)-(a_i*g_i,b_i,b_i,c_i)-(a_i,g_i,b_i,b_i)+(a_i*c_i,g_i*c_i,b_i*c_i,b_i*c_i).$$

Putting $g_i = (b_i - Ta_i)/(1 - T)$, we notice $a_i * g_i = b_i$. Hence, for any $i$, we may assume $a_i = b_i$. The condition $\partial_3(\sigma) = 0$ is thus formulated by

$$0 = \Big(\sum_i \alpha_i(a_i,a_i,a_i) - \alpha_i(a_i*c_i,a_i*c_i,a_i*c_i)\Big) + \sum_j \beta_j(d_j,f_j,f_j) - \beta_j(d_j*e_j,f_j,f_j). \tag{19}$$

One deals with the letter term $\sum_j \beta_j(d_j,e_j,f_j,f_j)$. We fix $\phi \in X$. Therefore $C^R_1(X;\mathbb{Z}) \ni 0 = \sum_j \beta_j(d_j - d_j*e_j) = \partial_2\big(\sum_j \beta_j(d_j,e_j)\big)$, where $j$ runs satisfying $f_j = \phi$. Since $H^R_2(X;\mathbb{Z}) \cong \mathbb{Z}$, the sum $\sum_j \beta_j(d_j,e_j)$ is homologous to $\gamma_\phi(\phi,\phi)$ for some $\gamma_\phi \in \mathbb{Z}$. Notice

$$\partial_5(a_j,b_j,c_j,\phi,\phi) = (a_j,c_j,\phi,\phi) - (a_j*b_j,c_j,\phi,\phi) - (a_j,b_j,\phi,\phi) + (a_j*c_j,b_j*c_j,\phi,\phi),$$

that is, $\partial_5(a_j, b_j, c_j, \phi, \phi) = \partial_3(a_j, b_j, c_j) \otimes (\phi, \phi)$. We therefore conclude that $\sum_j \beta_j(d_j, e_j, f_j, f_j)$ is homologous to $\sum_{\phi \in X} \gamma_\phi(\phi, \phi, \phi, \phi)$.

On the other hand, let us discuss the former term $\sum_i \alpha_i(a_i, a_i, a_i, c_i)$. By (19), we have $C_1^R(X; \mathbb{Z}) \ni 0 = \sum_i \alpha_i(a_i - a_i * c_i) = \partial_2(\sum_i \alpha_i(a_i, c_i))$. Since $H_2^R(X; \mathbb{Z}) \cong \mathbb{Z}$, the sum $\sum_i \alpha_i(a_i, c_i)$ is homologous to $\alpha(y, y)$ for some $\alpha \in \mathbb{Z}$ and $y \in X$. By definition, we note

$$\partial_5(a, a, a, d, f) = (a, a, a, f) - (a*d, a*d, a*d, f) - (a, a, a, f) + (a*f, a*f, a*f, d*f),$$

$$\partial_3(a, d, f) = (a, b, f) - (a*d, f) - (a, f) + (a*f, d*f).$$

In conclusion, $\sum_i \alpha_i(a_i, a_i, a_i, c_i)$ is homologous to $\alpha(y, y, y, y)$ as desired. $\square$

Finally, we will prove Lemma 5.7 as follows. For this, it is convenient to change another "coordinate system" of the complex $C_n^R(X; \mathbb{Z}[\mathrm{Inn}(X)])$ such as [M2, §2.1.3]. We denote elements of $\mathbb{Z}/e\mathbb{Z}$ by $\varepsilon$. We then define $C_n^{R_U}(X)$ to be the free $\mathbb{Z}$-module generated by elements $(\varepsilon, U_0; U_1, \ldots, U_n)$ of $\mathbb{Z}/e\mathbb{Z} \times X \times X^n$, and the boundary map to be

$$\partial_n(\varepsilon, U_0; U_1, \ldots, U_n) = \sum_{0 \le i \le n-1} (-1)^i \big((\varepsilon+1, TU_0; TU_1, \ldots, TU_{i-1}, TU_i + U_{i+1}, U_{i+2}, \ldots, U_n)$$
$$- (\varepsilon, U_0; U_1, \ldots, U_{i-1}, U_i + U_{i+1}, U_{i+2}, \ldots, U_n)\big). \quad (20)$$

We can see $\partial_{n-1} \circ \partial_n = 0$. Next, let us consider a bijection given by

$$\mathbb{Z}/e\mathbb{Z} \times X \times X^n \ni (\varepsilon, x_0; x_1, \ldots, x_n) \mapsto (\varepsilon, x_0 - x_1; x_1 - x_2, \ldots, x_{n-1} - x_n, x_n) \in \mathbb{Z}/e\mathbb{Z} \times X \times X^n.$$

The bijection induces a chain isomorphism from $C_n^R(X; \mathbb{Z}[\mathrm{Inn}(X)])$ to $C_n^{R_U}(X)$, where we use the identification $\mathrm{Inn}(X) \cong \mathbb{Z}/e\mathbb{Z} \ltimes X$ by Lemma 5.6.

*Proof of Lemma 5.7.* The proof is analogous to [N2, Theorem 6.1]. We put a chain homomorphism $\mathcal{Z} : C_n^{R_U}(X) \to C_n^{R_U}(X)$ defined by

$$\mathcal{Z}(\varepsilon, U_0; U_1, \ldots, U_n) := |X| \cdot (\varepsilon, 0; 0, \ldots, 0).$$

Furthermore, we define homomorphisms $D_{n,0}^j : C_n^{R_U}(X) \to C_{n+1}^{R_U}(X)$ and $D_{n,+}^j : C_n^{R_U}(X) \to C_{n+1}^{R_U}(X)$ by

$$D_{n,0}^j(\varepsilon, U_0; U_1, \ldots, U_n) = \sum_{y \in X} (\varepsilon, 0; 0, \ldots, 0, y, U_j - y, U_{j+1}, \ldots, U_n) \qquad \text{for } 0 \le j \le n,$$

$$D_{n,+}^j(\varepsilon, U_0; U_1, \ldots, U_n) = \begin{cases} \sum_{y \in X} (\varepsilon, 0; 0, \ldots, 0, y, -y, U_{j+1}, \ldots, U_n) & \text{for } 0 \le j < n, \\ \sum_{y \in X} (\varepsilon, 0; 0, \ldots, 0, y, -y) & \text{for } j = n. \end{cases}$$

In addition, we set $D_{n,0}^{n+1} = D_{n,+}^{n+1} = 0$. Then, by direct calculation, we can verify

$$\sum_{0 \le j \le n} (-1)^j \big(\partial_{n+1}(D_{n,0}^j + D_{n,+}^j) + (D_{n-1,0}^j + D_{n-1,+}^j)\partial_n\big) = (-1)^n \big(|X| \cdot \mathrm{id}_{C_n^{R_U}(X)} - \mathcal{Z}\big).$$

Then, we have the induced map $|X| \cdot \mathrm{id}_{H_n^{R_U}(X)} = (\mathcal{Z})_* : H_n^{R_U}(X) \to H_n^{R_U}(X)$. Therefore, for the proof, it suffices to show that any cycle of the form $\sum_\varepsilon a_\varepsilon (\varepsilon, 0; 0, \ldots, 0)$ is contained in the free subgroup of $C_n^{R_U}(X)$, where $a_\varepsilon \in \mathbb{Z}$. Indeed, the augmentation map $C_n^{R_U}(X) \to \mathbb{Z}$ is an $n$-cocycle and sends $(\varepsilon, 0; 0, \ldots, 0)$ to 1. □

# A  Appendix: $\pi_2(B^Q X) \otimes \mathbb{Z}_p$ of Alexander quandles $X$ on $\mathbb{F}_q$

This appendix computes $\pi_2(B^Q X) \otimes \mathbb{Z}_p$ for Alexander quandles of the forms $X = \mathbb{F}_q[T](T - \omega)$ (see §A.2). For this, in §A.1, we review some results in [M1, M2] and observe $b_i^Q = \dim_{\mathbb{F}_p}\bigl(H_i^Q(X;\mathbb{Z}) \otimes \mathbb{Z}_p\bigr)$.

## A.1  Review and some remarks on Mochizuki's cocycles

Mochizuki explicitly determined the second and third quandle cohomology groups [M1, M2]. To see this, we first recall from [M2, Theorem 2.2] that

$$\dim_{\mathbb{F}_q}(H_Q^2(X;\mathbb{F}_q)) = \#\bigl\{(i,j) \in \mathbb{Z}^2 \mid 1 \leq p^i < p^j < q,\ \omega^{p^i + p^j} = 1\bigr\} = \sum_{0 < i < h:\ \omega^{p^i+1}=1} h - i, \quad (21)$$

where $q = p^h$. In particular, one notices that $H_Q^2(X;\mathbb{F}_q)$ vanishes if and only if $\omega^{p^i+1} \neq 1$ for any $i < h$, that is, the order of $\omega$ is not divisible by $p^i + 1$ for any $i < h$.

Next, we discuss the third quandle cohomology $H_Q^3(X;\mathbb{F}_q)$, which is computed in [M2]. However, the statement of the main theorem [M2, Theorem 2.11] contained some minor mistakes. Then we will state the corrected statement. Recall the following polynomials over $\mathbb{F}_q$ in [M2, §2.2]:

$$F(a,b,c) := U_0^a \cdot U_1^b \cdot U_2^c, \qquad \chi(U_j, U_{j+1}) := \sum_{1 \leq i \leq p-1} (-1)^{i-1} i^{-1} U_j^{p-i} U_{j+1}^i,$$

$$E_0(a \cdot p, b) := \bigl(\chi(\omega U_0, U_1) - \chi(U_0, U_1)\bigr)^a \cdot U_2^b, \qquad E_1(a, b \cdot p) := U_0^a \cdot \bigl(\chi(U_1, U_2) - \chi(\omega^{-1} \cdot U_1, U_2)\bigr)^b.$$

Here for $(x_1, x_2, x_3) \in X^3$, we put $U_0 = x_1 - x_2$, $U_1 = x_2 - x_3$ and $U_2 = x_3$, similar to (20). Then these polynomials are regarded as functions of $C_3^Q(X;\mathbb{F}_q)$. Further, we let $\mathcal{Q}(q)$ be the set of quadruples $(q_1, q_2, q_3, q_4)$ satisfying the following conditions:

- $q_2 \leq q_3$, $q_1 < q_3$, $q_2 < q_4$, and $\omega^{q_1+q_3} = \omega^{q_2+q_4} = 1$. Here if $p = 2$, we omit $q_2 = q_3$.

- One of the following holds: **Case 1** $\omega^{q_1+q_2} = 1$

  **Case 2** $\omega^{q_1+q_2} \neq 1$, and $q_3 > q_4$.

  **Case 3** $(p \neq 2)$, $\omega^{q_1+q_2} \neq 1$, and $q_3 = q_4$.

  **Case 4** $(p \neq 2)$, $\omega^{q_1+q_2} \neq 1$, $q_2 \leq q_1 < q_3 < q_4$, $\omega^{q_1} = \omega^{q_2}$.

  **Case 5** $(p = 2)$, $\omega^{q_1+q_2} \neq 1$, $q_2 < q_1 < q_3 < q_4$, $\omega^{q_1} = \omega^{q_2}$.

For a quadruple $(q_1, q_2, q_3, q_4) \in \mathcal{Q}(q)$ satisfying the case 1, Mochizuki introduced $\Gamma(q_1, q_2, q_3, q_4) := F(q_1, q_2 + q_3, q_4) = U_1^{q_1} U_2^{q_2+q_3} U_3^{q_4}$. In other cases of $(q_1, q_2, q_3, q_4) \in \mathcal{Q}(q)$, see [M2, §2.2.2] for the definitions of polynomials $\Gamma(q_1, q_2, q_3, q_4)$. We state [M2, Theorem 2.11] as

**Theorem A.1** ([M2]). *The third quandle cohomology $H_Q^3(X; \mathbb{F}_q)$ is spanned by the following set $I_{q,\omega}$ composed of non-trivial 3-cocycles. Here $q_i$ means a power of the prime $p$ with $q_i < q$.*

$$\begin{aligned}
I_{q,\omega} := &\{F(q_1, q_2, q_3) \mid \omega^{q_1+q_2+q_3} = 1, \ q_1 < q_2 < q_3\} \cup \{F(q_1, q_2, 0) \mid \omega^{q_1+q_2} = 1, \ q_1 < q_2\} \\
&\cup \{E_0(p \cdot q_1, q_2) \mid \omega^{p \cdot q_1 + q_2} = 1, \ q_1 < q_2\} \cup \{E_1(q_1, q_2 \cdot p) \mid \omega^{q_1 + p \cdot q_2} = 1, \ q_1 \leq q_2\} \\
&\cup \{\Gamma(q_1, q_2, q_3, q_4) \mid (q_1, q_2, q_3, q_4) \in \mathcal{Q}(q)\}.
\end{aligned} \quad (22)$$

**Remark A.2.** In the original statement, the set $I_{q,\omega}$ contained certain 3-cocycles denoted by "$\Psi(a, q_1)$". However, we easily see that these 3-cocycles $\Psi(a, q_1)$ are cohomologous to zero (see [M2, §2.2.1]). This minor error results from the mistake "$B_{d-p^s}^{2(s)}(q)$" in the proof of [M2, Lemma 3.16] instead of "$B_{d-p^s}^2(q)$".

**Remark A.3.** By the universal coefficient theorem, we have $|I_{q,\omega}| = b_2^Q + b_3^Q$.

## A.2 Computations of $\pi_2(B^Q X) \otimes \mathbb{Z}_p$

We now compute $\pi_2(B^Q X) \otimes \mathbb{Z}_p$ in the cases of $\omega = -1$, of $q = p^{2h+1}$, and of $q = p^2$. In respect to a regular Alexander quandle $X$ of odd order, the following formula is useful:

$$\dim_{\mathbb{F}_p}(\pi_2(B^Q X) \otimes \mathbb{Z}_p) = b_3^Q - \frac{b_2^Q(b_2^Q - 1)}{2}, \quad (23)$$

where $b_n^Q$ is the dimension $\dim_{\mathbb{F}_p}(H_n^Q(X; \mathbb{Z}) \otimes \mathbb{Z}_p)$ with a prime $p > 2$. This formula is immediately obtained from Theorem 3.9.

### A.2.1 The case of $\omega = -1$

When $\omega = -1$ and $p$ is odd, let us compute $\pi_2(B^Q X) \otimes \mathbb{Z}_p$.

**Corollary A.4.** *Let $X$ be the product $h$-copies of the dihedral quandle $D_p$. Then*

$$\dim_{\mathbb{F}_p}(\pi_2(B^Q X) \otimes \mathbb{Z}_p) = \frac{h^2(h^2 + 11)}{12}.$$

*Proof.* We can verify that the order of the set $\mathcal{Q}(q)$ defined in §A.1 is $|\mathcal{Q}(q)| = h(h-1)(h+1)(5h-6)/24$. Hence, Theorem A.1 means $|I_{q,w}| = h(5h^3 - 6h^2 + 31h - 6)/24$. By (21), notice $b_2^Q = h(h-1)/2$. By (23), we thus conclude $\dim_{\mathbb{F}_p}(\pi_2(B^Q X) \otimes \mathbb{F}_p) = h^2(h^2 + 11)/12$. □

**Remark A.5.** We notice that the dimension is a quadruple function with respect to $h$. In particular, the second homotopy groups do not preserve the direct products of quandles.

### A.2.2 $\omega \neq -1$ and field extensions $\mathbb{F}_q$ of odd degrees

We next calculate $\dim(\pi_2(BX) \otimes \mathbb{F}_p)$, when $X = \mathbb{F}_q$ is an extension of odd degree and $\omega \neq -1$. To see this, we first show

**Lemma A.6.** *Let $\omega \neq \pm 1$ and $q = p^{2h+1}$ (possibly $p = 2$). Then $H_2^Q(X; \mathbb{Z}) \cong 0$.*

*Proof.* We first show $H_Q^2(X; \mathbb{F}_q) \cong 0$. For this, it suffices to show $\omega^{p^i+1} \neq 1$ for any $i \leq 2h$ by (21). A key is that, if $p > 2$, then $(p^{2m}+1)/2$ and $(p^{2n-1}+1)/2$ are relatively prime for any $n, m \in \mathbb{Z}$; This is easily verified by induction on $n+m$. We here assume $\omega^{p^i+1} = 1$ for some $i \leq 2h$. Then $\omega^{p^{2h-i+1}+1} = (\omega^{p^i+1} \omega^{p^{2h+1}-1})^{-p^i} = 1$. This means that $p^i + 1$ and $p^{2h-i+1} + 1$ have a common divisor except for 2, which contradicts the key. A similar discussion holds for the case $p = 2$. Hence $H_Q^2(X; \mathbb{F}_q) \cong 0$. Finally, since $H_2^Q(X; \mathbb{Z})$ is annihilated by $q$ (see [N2, Corollary 6.2]), we have $H_2^Q(X; \mathbb{Z}) \cong 0$ by the universal coefficient theorem. □

Therefore, combing Lemma A.6 with Theorem A.1 immediately concludes

**Corollary A.7.** *If $\mathbb{F}_q$ is of odd degree and $\omega \neq -1$, then $\dim(\pi_2(BX) \otimes \mathbb{F}_p)$ is equal to the order $|\{(q_1, q_2, q_3) \mid \omega^{q_1+q_2+q_3} = 1\}|$.*

**Example A.8.** Let $q = p$. If $\omega \neq -1$, then the second and third quandle cohomologies vanish; hence $\pi_2(B^Q X) \cong 0$.

**Example A.9.** For example, we let $q = p^3$. One notice that, when $\omega^{1+p+p^2} \neq 1$, $H_Q^3(X; \mathbb{F}_q)$ vanishes; hence, $\pi_2(B^Q X) \cong 0$. On the other hand, if $\omega^{1+p+p^2} = 1$, then $I_{q,\omega}$ consists of the polynomial $F(1, p, p^2)$. Hence, $\pi_2(B^Q X) \otimes \mathbb{Z}_p \cong \mathbb{Z}_p$

We further consider an example of $(p, w) = (7, 2)$, which satisfies $\omega^{1+p+p^2} = 1$. Let $X = (\mathbb{Z}[T]/(p, T - \omega))^3$. Hence $\pi_2(B^Q X) \otimes \mathbb{Z}_7 \cong \mathbb{Z}_7$. Compare this with Example A.8. Namely, direct products of quandles sometimes produce non-trivial homotopy groups of $B^Q X$. There are many examples of such pairs $(p, w)$, e.g., (13,3), (19,7), (37,10) and so on.

### A.2.3 The case of $q = p^2$ with $\omega \neq -1$

We here compute $\pi_2(B^Q X)$ of extension fields $X$ over $\mathbb{F}_p$ of degree 2. Without $\omega = -1$, the dimensions of $\pi_2(B^Q X) \otimes \mathbb{Z}_p$ are determined by the following three cases.

**Example A.10.** We consider the case of $\omega^{1+p} \neq 1$. Then we have $b_2^Q = 0$ by (21). Further, by Theorem A.1, we obtain $b_3^Q = 0$. Therefore, $\pi_2(B^Q X) \cong 0$.

**Example A.11.** We consider the case where $\omega \neq -1$, $\omega^{1+p} = 1$ and $p \neq 2$. Then the equalities (21) mean $b_2^Q = 1$. Further, by Theorem A.1, the third cohomology $H_Q^3(X; \mathbb{F}_q)$ is spanned by

$$I_{p^2, \omega} = \{F(1, p, 0),\ E_1(1, p),\ E_1(p, p^2),\ \Gamma(1, p, 1, p)\}.$$

Hence $b_3^Q = 3$. Therefore we have $\dim_{\mathbb{F}_p}(\pi_2(B^Q X) \otimes \mathbb{Z}_p) = 3$ by the formula (23).

More generally, for given a quandle of the form $X = \mathbb{F}_q[T]/(T - \omega)$ with $q > p^2$, we can compute $\dim(\pi_2(B^Q X) \otimes \mathbb{Z}_p)$ in a similar manner.

**Acknowledgment**

The author is grateful to Takuro Mochizuki for helpful communications on quandle cocycles and Remark A.2. He expresses his gratitude to Syunji Moriya for valuable conversations on rational homotopy theory. He is grateful to J. Scott Carter, Inasa Nakamura, Masahico Saito and Kokoro Tanaka for useful discussions on Remark 3.3 and the paper [CEGS]. He also thanks Tomotada Ohtsuki and Dai Tamaki for valuable comments.

Research Institute for Mathematical Sciences, Kyoto University, Sakyo-ku, Kyoto, 606-8502, Japan

E-mail address: `nosaka@kurims.kyoto-u.ac.jp`